%% file: BCM18.tex
\documentclass[10pt,final]{siamltex}
\pdfoutput=1
\usepackage{xcolor}
\usepackage{amsmath}
\usepackage{amssymb}
\usepackage{mathrsfs}
\usepackage{graphicx}
\usepackage{bm}
\usepackage{ucs}
\usepackage[utf8x]{inputenc}
\usepackage{wrapfig}
\usepackage{color}
\usepackage{float}
\usepackage{epstopdf}

\usepackage[breaklinks,colorlinks=true,linkcolor=blue,citecolor=red,
  backref=page]{hyperref}

\newtheorem{remark}{Remark}
\newtheorem{algorithm}[theorem]{Algorithm}
\input{mathmacros.tex}

\begin{document}

%%%% Article title to be placed here

\title{A direct approach to imaging in a waveguide with perturbed geometry}
\author{Liliana Borcea\footnotemark[1],  Fioralba
  Cakoni\footnotemark[2] and Shixu Meng\footnotemark[3]}

\maketitle

\renewcommand{\thefootnote}{\fnsymbol{footnote}}
\footnotetext[1]{Department of Mathematics, University of Michigan,
  Ann Arbor, MI 48109. {\tt borcea@umich.edu}}
\footnotetext[2]{Department of Mathematics, Rutgers University, New
  Brunswick, NJ 08901. {\tt fc292@math.rutgers.edu}}
\footnotetext[3]{Department of Mathematics, University of Michigan,
  Ann Arbor, MI 48109. {\tt shixumen@umich.edu}}

\begin{abstract}
We introduce a direct, linear sampling approach to imaging in an
acoustic waveguide with sound hard walls.  The waveguide terminates at
one end and has unknown geometry due to compactly supported wall
deformations. The goal of imaging is to determine these deformations
and to identify localized scatterers in the waveguide, using a remote
array of sensors that emits time harmonic probing waves and records
the echoes. We present a theoretical analysis of the imaging approach and
illustrate its performance with numerical simulations.
\end{abstract}
\begin{keywords}
Linear sampling method, waveguide, inverse scattering
\end{keywords}

\section{Introduction and formulation of the problem}
\label{sect:Intro}
Sensor array imaging in waveguides has applications in underwater
acoustics \cite{Wirgi2000,baggeroer1993overview}, nondestructive
evaluation of slender structures
\cite{chillara2015review,rizzo2010ultrasonic}, imaging of and in
tunnels \cite{schultz2007remote,haack1995state,bedford2014modeling},
etc. It is a particular inverse wave scattering problem that has been
studied extensively for waveguides with known and simple geometry.
The wave equation in such empty waveguides can be solved with
separation of variables and the wave field is a superposition of
propagating, evanescent and possibly radiating modes that do not
interact with each other.  A sample of the existing mathematical
literature is
\cite{dediu2006recovering,bourgeois2008linear,monk2012sampling,
  monk2016inverse,borcea2015imaging,tsogka2013selective,tsogka2017imaging}
and examples of imaging with experimental validation are in
\cite{mordant1999highly,philippe2008characterization}.

The problem is more difficult when the waveguide has variable and
unknown geometry. Studies of wave propagation in waveguides with
random boundary \cite{alonso2011wave,
  borcea2014paraxial,Gomez1,borcea2017transport,borcea2017pulse} show
that even small amplitude fluctuations of the walls can have a
significant scattering effect (i.e., mode coupling) over long
distances of propagation, manifested by the randomization of the wave
field.  While experiments like time reversal
\cite{garnier2007pulse,borcea2014paraxial} take advantage of such net
scattering, the uncertainty of the boundary poses a serious impediment
to imaging that has lead to proposals of new data processing and
measurement setups \cite{borcea2010source,Gomez,
  borcea2014paraxial,acosta2015source,borcea2018ghost}.

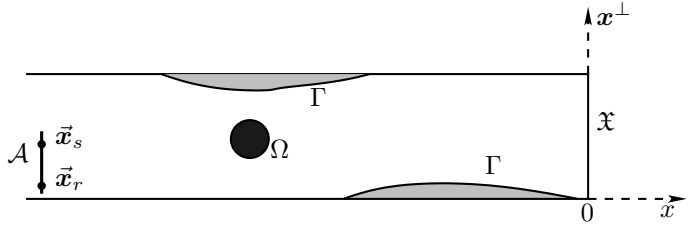
\begin{figure}[t]
\begin{center}
\includegraphics[width=0.68\textwidth]{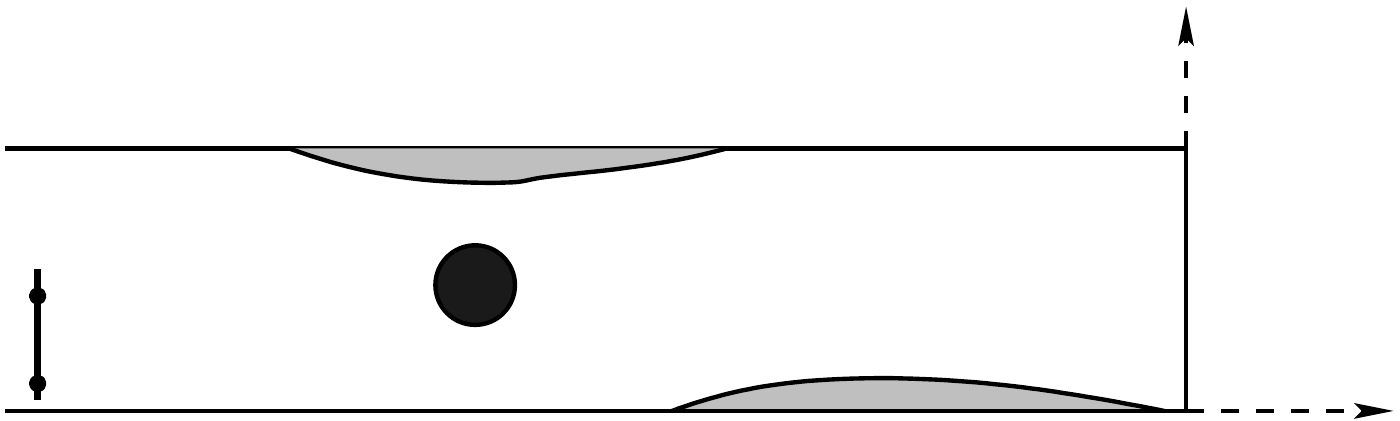}%
\end{center}
\setlength{\unitlength}{3947sp}%
\begingroup\makeatletter\ifx\SetFigFont\undefined%
\gdef\SetFigFont#1#2#3#4#5{%
  \reset@font\fontsize{#1}{#2pt}%
  \fontfamily{#3}\fontseries{#4}\fontshape{#5}%
  \selectfont}%
\fi\endgroup%
\begin{picture}(9591,2871)(1918,-3049)
\put(2800,100){\makebox(0,0)[lb]{\smash{{\SetFigFont{7}{8.4}{\familydefault}{\mddefault}{\updefault}{\color[rgb]{0,0,0}{\normalsize $\mathcal{A}$ }}%
}}}}
\put(3100,-50){\makebox(0,0)[lb]{\smash{{\SetFigFont{7}{8.4}{\familydefault}{\mddefault}{\updefault}{\color[rgb]{0,0,0}{\normalsize $\vx_r$ }}%
}}}}
\put(3100,210){\makebox(0,0)[lb]{\smash{{\SetFigFont{7}{8.4}{\familydefault}{\mddefault}{\updefault}{\color[rgb]{0,0,0}{\normalsize $\vx_s$ }}%
}}}}
\put(6900,-250){\makebox(0,0)[lb]{\smash{{\SetFigFont{7}{8.4}{\familydefault}{\mddefault}{\updefault}{\color[rgb]{0,0,0}{\normalsize $x$ }}%
}}}}
\put(6500,960){\makebox(0,0)[lb]{\smash{{\SetFigFont{7}{8.4}{\familydefault}{\mddefault}{\updefault}{\color[rgb]{0,0,0}{\normalsize $\bxp$ }}%
}}}}
\put(6400,-270){\makebox(0,0)[lb]{\smash{{\SetFigFont{7}{8.4}{\familydefault}{\mddefault}{\updefault}{\color[rgb]{0,0,0}{\normalsize $0$ }}%
}}}}
\put(4700,450){\makebox(0,0)[lb]{\smash{{\SetFigFont{7}{8.4}{\familydefault}{\mddefault}{\updefault}{\color[rgb]{0,0,0}{\normalsize $\GD$ }}%
}}}}
\put(4450,120){\makebox(0,0)[lb]{\smash{{\SetFigFont{7}{8.4}{\familydefault}{\mddefault}{\updefault}{\color[rgb]{0,0,0}{\normalsize $\Omega$ }}%
}}}}
\put(6500,300){\makebox(0,0)[lb]{\smash{{\SetFigFont{7}{8.4}{\familydefault}{\mddefault}{\updefault}{\color[rgb]{0,0,0}{\normalsize $\cX$ }}%
}}}}
\put(5800,20){\makebox(0,0)[lb]{\smash{{\SetFigFont{7}{8.4}{\familydefault}{\mddefault}{\updefault}{\color[rgb]{0,0,0}{\normalsize $\GD$ }}%
}}}}
\end{picture}%
\vspace{-2.3in}
\caption{Illustration of the imaging setup in a terminating
  waveguide. The system of coordinates is $\vx = (x,\bxp)$ with range
  $x$ measured from the end wall and cross-range $\bxp$ in the
  cross-section $\cX$ of the waveguide. The wall deformation of the
  waveguide is modeled by the boundary $\GD$ of the domain $\cD$
  drawn in gray.  A localized scatterer supported in $\Omega$ is drawn
  in black. The array of sensors lies in the set $\cA$. The source
  and receiver locations are denoted by $\vx_s$ and 
  $\vx_r$.}
\label{fig:setup}
\end{figure}

Here we consider a different type of wall deformations, with larger
amplitude but compact support, and pursue a linear-sampling approach
for estimating these deformations and localized scatterers in the
waveguide.  Motivated by the application of imaging in tunnels, we
consider a waveguide that terminates, as illustrated in Figure
\ref{fig:setup}. For simplicity, we limit the study to acoustic waves
and to sound hard walls, but the linear sampling approach can be
extended to other boundary conditions and to electromagnetic and
elastic waves. We refer to  \cite{bourgeois2011use,bourgeois2013use,yang2015scattering} for linear sampling
imaging in waveguides with elastic waves  and to \cite{yang2015scattering} for imaging with electromagnetic waves.

Let us denote by $\cW_o$ the ideal waveguide with unperturbed walls
modeled by the boundary $\partial \cW_o$, and use the system of
coordinates $\vx = (x,\bx^\perp) \in \RR^{d}$ shown in Figure
\ref{fig:setup}, with range $x$ measured along the axis of $\cW_o$,
starting from the end wall.  The cross-range coordinates $\bxp$ lie in
the cross-section of $\cW_o$, denoted by $\cX \subset
\RR^{d-1}$. This is a compact Lipschitz domain when $d = 3$, or an
interval of finite length $|\cX|$ when $d = 2$. In our system of
coordinates we have
\begin{equation}
\cW_o = (-\infty, 0) \times \cX, \qquad \partial \cW_o =
\Big((-\infty,0)\times \partial \cX\Big) \cup \Big(\{0\} \times \cX\Big),
\label{eq:Int1}
\end{equation}
and we model the unknown waveguide by 
\begin{equation}
\cW = \cW_o \cap (\mathbb{R}^d \setminus \overline{\cD}),
\label{eq:Int2}
\end{equation}
where $\cD$ is a Lipschitz domain compactly supported in the sector
$(-x_\star,0) \times \cX$ of $\cW_o$, with part of the boundary
$\partial \cD$ lying in $\partial \cW_o$. We denote this part by
$\Gamma_o$ and model the unknown waveguide walls by 
\begin{equation}
\GD = \partial \cD \setminus 
\overline{\Gamma}_o \subset \cW_o,
\label{eq:Int3}
\end{equation} 
where the bar denotes the closure of $\Gamma_o$. The waveguide is
filled with a homogeneous medium (e.g. air) but it may contain one or
more impenetrable or penetrable scatterers supported in the compact
set $\Omega$, satisfying
\begin{equation}
\Omega  \subset \cW \cap \Big((-x_\star,0) \times \cX\Big).
\label{eq:Int4}
\end{equation} 
This is  a Lipschitz domain or the union of a few
disjoint such domains.

The imaging problem is to estimate $\GD$ and $\Omega$ using data
gathered by an array of $\JA$ sensors located in the set 
\begin{equation}
\cA
\subseteq \{\xA\} \times \cX, \quad \xA < x_\star < 0,
\label{eq:Int5}
\end{equation} 
called the array aperture.  The array
probes the waveguide by emitting a time harmonic wave from one of the
sensors, at location $\vx_s$, and measures the echoes $u(\vx_r,\vx_s)$
at all the sensors $\{\vx_r\}_{r = 1,\ldots, J_{_\cA}}$. Although $s$
and $r$ are indexes in the set $\{1, \ldots, \JA\}$, we use them
consistently to distinguish between the source and receiver. The data
gathered successively, with one source at a time, form the $\JA \times
\JA$ response matrix $\big(u(\vx_r;\vx_s)\big)_{1 \le r,s\le
  \JA}$. The goal is to show with analysis and numerical
simulations how  the linear sampling approach
estimates $\Gamma$ and $\Omega$ from this matrix.

The paper is organized as follows: We begin in  {Section}
\ref{sect:Problem1} with the estimation of $\GD$. The estimation of
both $\GD$ and $\Omega$ is considered in  {Section} \ref{sect:Problem2}. The assessment with numerical simulations is in
 {Section} \ref{sect:Numerics}. We end with a summary in section
\ref{sect:Summary}.

\section{Imaging wall deformations}
\label{sect:Problem1}
We define in  {Section} \ref{sect:Prob1.1} the Green's function in the
unperturbed waveguide, which models the incident wave emitted by a
source in the array. The model of the scattered wave measured at the
array is given in  {Section} \ref{sect:Prob1.0}.  The linear sampling
approach is  analyzed in {Section} \ref{sect:Prob1.2}, for
the case of a full aperture array. Imaging with a partial aperture
array is described in {Section} \ref{sect:Prob1.3}.

\subsection{The incident wave field}
\label{sect:Prob1.1}
Let us denote by $G(\vx,\vy)$ the Green's function in the ideal
waveguide $\cW_o$, for an arbitrary source location $\vy = (y,\byp) \in
\cW_o$.  The model of the incident wave emitted by the source at
location $\vx_s \in \cA$ is then
\begin{equation}
u^{{\rm inc}}(\vx,\vx_s) = G(\vx,\vx_s).
\label{eq:1.0}
\end{equation}

The Green's function satisfies the Helmholtz equation
\begin{equation}
\big(\Delta_{\vx} + k^2 \big) G(\vx,\vy) = - \delta(\vx - \vy),
\quad \vx \in \cW_o, \label{eq:1.1}
\end{equation}
where $\Delta_{\vx}$ is the Laplacian with respect to $\vx$ and $k$ is
the wavenumber.  At the sound hard walls $\partial \cW_o$ we have the
boundary condition
\begin{equation}
 \frac{\partial G(\vx,\vy)}{\partial \vnu}
= 0, \quad \vx \in \partial W_o, \label{eq:1.2}
\end{equation}
where $\vnu$ denotes the outer unit normal at $\vx$, and for $\vx \in
\cW_o$ with range coordinate $x < y$ we impose the radiation
condition  {formulated precisely in Definition \ref{def.1}}, which states that $G(\vx,\vy)$ is a bounded and outgoing
wave.

Due to the simple geometry of $\cW_o$, the Green's function can be
written explicitly using the eigenfunctions $\{\psi_j(\bx^\perp)\}_{j \ge
  0}$ of the Laplacian $\Delta_{\bx^\perp}$ in $\cX$, satisfying
\begin{align}
-\Delta_{\bxp} \psi_j(\bxp) &= \lambda_j \psi_j(\bxp), \qquad \bxp \in
\cX, \nonumber\\ \frac{\partial \psi_j(\bxp)}{\partial \nu_{{\bxp}}}
&= 0, \quad \bxp \in \partial \cX,
\label{eq:1.3}
\end{align}
where $\nu_{{\bxp}}$ is the outer normal at $\bxp$, in the plane of
$\cX \subset \RR^{d-1}$.  The spectral theorem for compact
self-adjoint linear operators \cite[Theorem 2.36]{mclean2000strongly}
implies that these eigenfunctions form a complete orthonormal basis of
$L^2(\cX)$ and that the eigenvalues $\lambda_j$ are real and
non-negative.  The first eigenvalue $\la_o = 0$ is simple and
corresponds to the constant eigenfunction $\psi_0(\bxp) ={1}/\sqrt{|\cX|}$.
The other eigenvalues satisfy
\begin{equation}
0 = \lambda_o < \lambda_1 \le \lambda_2 \le \ldots, \qquad \lim_{j \to
  \infty} \lambda_j = \infty.
\label{eq:1.4}
\end{equation}
The expression of the Green's function is 
\begin{equation}
G(\vx,\vy) = \sum_{j=0}^\infty \frac{i}{2 \beta_j} \psi_j(\byp) \psi_j(\bxp)
\Big( e^{i \beta_j |x-y| } + e^{i \beta_j|x+y|} \Big), 
\label{eq:1.5}
\end{equation}
where 
\begin{equation}
\beta_j = \left\{ \begin{array}{ll}
\sqrt{k^2 - \lambda_j}, \quad &j = 0, 1, \ldots, J, \\ \\
i \sqrt{\lambda_j - k^2}, & j > J, 
\end{array} \right. 
\label{eq:1.6}
\end{equation}
and $J$ is the largest index $j$ such that $\lambda_j \le  k^2$. 

Note that at points $\vx = (x,\bxp)\in \cW_o$ between the source at
$\vy = (y,\byp)$ and the end wall i.e., for range $x \in (y,0)$, the
expression \eqref{eq:1.5} consists of $J+1$ propagating modes
$\{\psi_j(\bxp) e^{\pm i \beta_j x}\}_{0 \le j \le J}$ and infinitely
many growing and decaying (evanescent) modes $\{\psi_j(\bxp) e^{\pm
  \beta_j x}\}_{j > J}$ with complex amplitudes that depend on
$\vy$. The propagating modes can be understood as superpositions of
plane waves with wave vector $(\pm \beta_j, \boldsymbol{\kappa}_j)$,
where $ \boldsymbol{\kappa}_j \in \RR^{d-1}$ has the square Euclidian 
norm $\lambda_j$. These waves propagate forward and backward in the
range direction, at group speed
\[
c\Big(\frac{d \beta_j}{d k}\Big)^{-1} = c \frac{\beta_j}{k}, \qquad j =
0, \ldots, J,
\]
where $c$ is the wave speed in the homogeneous medium that fills the
waveguide. The fastest mode indexed by $j = 0$ propagates at
speed $c$. The slowest mode corresponds to $j = J$ and we assume that
$\lambda_J < k^2$, so that $\beta_J \ne 0$. The wavenumber is
imaginary for indexes $j > J$ and the modes grow or decay
exponentially in range.

At points $\vx$ with range coordinate $x < y$, the expression
\eqref{eq:1.5} consists of $J+1$ outgoing (backward) propagating modes
$\{\psi_j(\bxp) e^{- i \beta_j x}\}_{0 \le j \le J}$ and infinitely
many decaying (evanescent) modes $\{\psi_j(\bxp) e^{\beta_j x}\}_{j> J}$. This is
the explicit statement of the radiation condition for the Green's
function.

\subsection{The array response matrix}
\label{sect:Prob1.0}
The scattered field $u(\vx,\vx_s)$ due to the incident wave
\eqref{eq:1.0} is the function in $H_{{\rm loc}}^1(\cW)$ satisfying
the Helmholtz equation
\begin{equation}
\big(\Delta_{\vx} + k^2 \big) u(\vx,\vx_s) = 0, \qquad \vx \in \cW,
\label{eq:1.7}
\end{equation}
with the Neumann boundary conditions 
\begin{align}
 \frac{\partial u(\vx,\vx_s)}{\partial \vnu} &= 0, \quad   \partial W_o \setminus \overline{\Gamma}_o, \\ \frac{\partial
   u(\vx,\vx_s)}{\partial \vnu} &= - \frac{\partial
   G(\vx,\vx_s)}{\partial \vnu}, \qquad \vx \in \Gamma,
\label{eq:1.7b}
\end{align}
at the sound hard walls, and the radiation condition at points $\vx
\in \cW$ with range coordinate $x < x_\star$. Due to the assumption
that the wall deformation is supported in the range interval
$(x_\star,0)$, with $\xA < x_\star$, the radiation condition is as in
the previous section:
\begin{definition}
\label{def.1}
The radiation
condition at points $\vx = (x,\bxp) \in \cW$ with $x < x_\star$ means
that $u(\vx,\vx_s)$ is a superposition of $J+1$ backward going modes and infinitely many
decaying modes,
\begin{equation}
u(\vx,\vx_s) = \sum_{j=0}^\infty \alpha_j(\vx_s,\Gamma) \psi_j(\bxp) e^{-i
  \beta_j x}, \qquad \vx = (x,\bxp), ~ ~ x < x_\star.
\label{eq:1.8}
\end{equation}
Each term (mode) in the sum is a special solution of the
Helmoltz equation in the sector $(-\infty,x_\star) \times \cX$ of
$\cW$.  The complex amplitudes $\alpha_j$ depend on
$\vx_s$ and $\Gamma$.
\end{definition}

The array is located far from the wall deformation, so the response
matrix can be modeled as 
\begin{equation}
u(\vx_r,\vx_s) \approx \sum_{j=0}^J \alpha_j(\vx_s,\Gamma) e^{-i
  \beta_j \xA} \psi_j(\bxp_r), \qquad \forall \, \vx_r,\vx_s \in \cA,
\label{eq:1.9}
\end{equation}
where we neglect the evanescent waves. 

\subsection{The linear sampling approach}
\label{sect:Prob1.2}
In this section we show how to use the linear sampling approach to
estimate $\GD$ from the array response matrix with entries
\eqref{eq:1.9}. In the analysis we assume that the sensors are located
very close together in the array and we replace sums over the sensor
indexes by integrals over $\cA$. Although we keep the notation $\vx_s$
and $\vx_r$ for the source and receiver locations, these are now
vectors that vary continuously in $\cA$.  We begin with the case of
full array aperture
\begin{equation}
\cA = \{\xA\} \times \cX,
\label{eq:FA}
\end{equation} 
and postpone until the next section the discussion for partial
aperture. {However we remark that  the theoretical justification  of the linear sampling method for partial aperture remains unchanged.}

\subsubsection{Analysis of the linear sampling approach}
Let us introduce the so-called near field integral operator $N :
L^2(\cA) \to L^2(\cA)$ defined by
\begin{equation}
N g(\vx_r) = \int_{\cA} dS_{\vx_s} \, u(\vx_r,\vx_s) g(\vx_s), \qquad
\forall g \in L^2(\cA), ~ ~ \vx_r \in \cA,
\label{eq:1.10}
\end{equation}
where we note that the assumption \eqref{eq:FA} implies that the
cross-range components of $\vx_r, \vx_s$ lie in $\cX$.  By linear
superposition, the function $Ng(\vx_r)$ represents the scattered wave
received at $\vx_r$, due to an illumination $g(\vx_s)$ from all the
source points $\vx_s \in \cA$. The linear sampling method uses this
$g(\vx_s)$ as a control at the array, which focuses the wave at a
point $\vz$ in the imaging domain, so that the received wave
$Ng(\vx_r)$ equals $G(\vx_r,\vz)$. It turns out that the control
function $g$ is not physical (i.e., it is not bounded in $L^2(\cA)$)
if $\vz \notin \cD$, and this leads to the linear sampling imaging
approach.

Our analysis of the linear sampling method is based on the following
factorization of the near field operator, proved in appendix
\ref{ap:A}:

\vspace{0.05in}
\begin{lemma}
\label{lem.1}
The operator $N$ defined in \eqref{eq:1.10} has the factorization
\begin{equation}
N = \TDA \TAD,
\label{eq:1.11}
\end{equation}
where $\TAD: L^2(\cA) \to H^{-\frac{1}{2}}(\GD)$ is the operator 
\begin{equation}
\TAD g (\vz) = \partial_{\vnuz} \hspace{-0.03in} \int_{\cA} d S_{\vx_s}\,
G(\vz,\vx_s) g(\vx_s), \qquad \forall g \in L^2(\cA), ~ \vz \in \GD,
\label{eq:1.12}
\end{equation}
and $\TDA:H^{-\frac{1}{2}}(\GD)\to L^2(\cA)$ is the operator defined by the
trace $\TDA f = w|_{\cA}$ of the solution of
\begin{align}
\big(\Delta_{\vx} + k^2 \big) w(\vx) &= 0, \qquad \vx \in
\cW, \label{eq:1.13a} \\ \frac{\partial w(\vx)}{\partial \vnu} &= 0,
\quad \vx \in \partial W_o \setminus \overline{\Gamma}_o, \\ \frac{\partial
  w(\vx)}{\partial \vnu} &= - f(\vx), \qquad \vx \in
\Gamma, \label{eq:1.13d}
\end{align}
satisfying a radiation condition as in Definition \ref{def.1}.
\end{lemma}

\vspace{0.05in} 
We conclude from the factorization \eqref{eq:1.11} that 
\begin{equation}
\mbox{range} \big(N\big) \subset \mbox{range} \big(\TDA\big) \subset
L^2(\cA).
\label{eq:1.14}
\end{equation}
We also see from \eqref{eq:1.13a}--\eqref{eq:1.13d} that the range of
$\TDA$ consists of traces {on $\cA$} of functions that satisfy Helmholtz's
equation in $\cW$ with homogeneous Neumann boundary conditions on
$ { \partial W_o \setminus \overline{\Gamma}_o}$ and the radiation condition. An example of such a
function is $G(\vx,\vz)$ for any $\vz \in \cD$. The next lemma, proved
in appendix \ref{ap:A}, uses this observation to distinguish between
points inside and outside $\cD$.

\vspace{0.05in} 
\begin{lemma}
\label{lem.2}
Let $\vz $ be a search point in $\cW_o$,
between the array and the end wall. Then, $\vz \in \cD$ if and only if
$G(\cdot,\vz)|_{\cA} \in \mbox{range} \big(\TDA\big)$.
\end{lemma}

\vspace{0.05in} Since $\TDA$ is unknown, we cannot determine the
support of $\cD$ directly from Lemma \ref{lem.2}. We only know the
near field operator \eqref{eq:1.10} with range satisfying
\eqref{eq:1.14}. While $G(\cdot,\vz)|_{\cA} \in \mbox{range}
\big(\TDA\big)$ implies the existence of $f \in H^{-\frac{1}{2}}(\GD)$ such
that $\TDA f = G(\cdot,\vz)|_{\cA}$, it is not clear that $f$ is in $\mbox{range} \big(\TAD\big)$. 
The next lemma, proved in appendix \ref{ap:A},
shows that  $f$ can be approximated arbitrarily well by some
$\widetilde f \in \mbox{range} \big(\TAD\big)$ and, furthermore, 
that $N \widetilde f \approx G(\cdot,\vz)|_{\cA}$.

\vspace{0.05in} 
\begin{lemma}
\label{lem.3}
The linear operator $\TAD$ is bounded and has dense range in
$H^{-\frac{1}{2}}(\GD)$.  The linear operator $\TDA$ is compact and has dense
range in $L^2(\cA)$.
\end{lemma}

\vspace{0.05in} Gathering the results in Lemmas \ref{lem.1}--\ref{lem.3}, 
we can now prove the following result for the linear sampling approach:

\vspace{0.05in} 
\begin{theorem}
\label{thm.1}
Let $ \vz$ be a search point in $\cW_o$,
between the array and the end wall. For any $\ep > 0$ {let $g_{\vz}^\ep \in L^2(\cA)$ satisfy}
\begin{equation}
\|Ng_{\vz}^\ep - G(\cdot,\vz)\|_{L^2(\cA)} < \ep.
\label{eq:1.15}
\end{equation}
{(which obviously exists since the range of $N$ is dense in $L^2(\cA)$)}. 

There are two possibilities:

\vspace{0.03in}
\begin{enumerate}
\item If $\vz \in \cD$, {there exists a} $g_{\vz}^\ep$ satisfying
  \eqref{eq:1.15}  {such that}  the norm $\|\TAD g_{\vz}^\ep\|_{H^{-\frac{1}{2}}(\GD)}$
  remains bounded as $\ep \to 0$.
\item If $\vz \notin \cD$, for any $g_{\vz}^\ep$ satisfying
  \eqref{eq:1.15}, $\displaystyle \lim_{\ep \to 0} \|\TAD
  g_{\vz}^\ep\|_{H^{-\frac{1}{2}}(\GD)} = \infty$.
\end{enumerate}
\end{theorem}

\vspace{0.05in} This theorem says that it is possible to estimate the
support of $\cD$ and therefore the deformed walls $\GD$, from the
magnitude of $\|\TAD g_{\vz}^\ep\|_{H^{-\frac{1}{2}}(\GD)}$. However, this
norm cannot be computed, because we do not know $\GD$ and therefore
$\TAD$. To obtain an imaging method, we use instead the norm
$\|g_{\vz}^\ep\|_{L^2(\cA)}.$ Recalling from Lemma \ref{lem.3} that
$\TAD$ is a bounded linear operator, we have
\begin{equation}
\|g_{\vz}^\ep\|_{L^2(\cA)} \ge \frac{\|\TAD
  g_{\vz}^\ep\|_{H^{-\frac{1}{2}}(\GD)}}{\|\TAD\|},
\label{eq:1.16}
\end{equation}
so if $z \notin \cD$, we conclude from case 2.  of Theorem \ref{thm.1}
that $\displaystyle\lim_{\ep \to 0}\|g_{\vz}^\ep\|_{L^2(\cA)} =
\infty$. However, if $z \in \cD$ we cannot guarantee that
$\|g_{\vz}^\ep\|_{L^2(\cA)}$ remains bounded, because there may be
large components of $g_{\vz}^\ep$ in the null space of
$\TAD$. Nevertheless, we can control such components by searching for
the minimum norm solution $g_{\vz}^\ep$ of \eqref{eq:1.15} or,
similarly, by minimizing $\|N g - G(\cdot,\vz)\|_{L^2(\cA)}$ using
Tikhonov regularization, as explained in section \ref{sect:ImagAlg}.

\textbf{Proof of Theorem \ref{thm.1}:} Let us begin with case 1., for
search point $\vz \in \cD$. By Lemma \ref{lem.2}, we conclude that
$\exists \, f_{\vz} \in H^{-\frac{1}{2}}(\GD)$ such that 
\begin{equation}
\TDA f_{\vz}(\vx) = G(\vx,\vz)|_{\cA}.
\label{eq:pf1}
\end{equation}
By Lemma \ref{lem.3}, since $\mbox{range}\big({\TAD}\big)$ is dense in
$H^{-\frac{1}{2}}(\GD)$, for any $\ep > 0$ there exists $ g_{\vz}^\ep \in
L^2(\cA)$ such that
\begin{equation}
\|\TAD g_{\vz}^\ep - f_{\vz}\|_{H^{-\frac{1}{2}}(\GD)} < \frac{\ep}{\|\TDA\|},
\label{eq:pf3}
\end{equation}
where we used that $\TDA$ is bounded, per Lemma \ref{lem.3}. Then, the
factorization in Lemma \ref{lem.1} and \eqref{eq:pf1} give that this
$g_{\vz}^\ep$ satisfies
\begin{align}
\|N g_{\vz}^\ep - G(\cdot,\vz) \|_{L^2(\cA)} = \big\|\TDA \big(\TAD
g_{\vz}^\ep - f_{\vz}\big) \|_{L^2(\cA)} < \ep.
\label{eq:pf4}
\end{align}
We also have using the triangle inequality in \eqref{eq:pf3} that 
\[
\|\TAD g_{\vz}^\ep \|_{H^{-\frac{1}{2}}(\GD)} \le \frac{\ep}{\|\TDA\|}
+ \|f_{\vz}\|_{H^{-\frac{1}{2}}(\GD)} \to
\|f_{\vz}\|_{H^{-\frac{1}{2}}(\GD)} ~ ~ \mbox{ as }~ \ep \to 0.
\]
This proves case 1. of the theorem.

For case 2., let $\vz \notin \cD$ and conclude from Lemmma \ref{lem.2}
that $\forall f \in H^{-\frac{1}{2}}(\GD)$,
\begin{equation}
\|\TDA f - G(\cdot,\vz)\|_{L^2(\cA)} > 0.
\label{eq:PF0}
\end{equation}
Nevertheless, since $G(\vx,\vz) \in L^2(\cA)$ for $\vz \notin \cA$ and
$\mbox{range}\big(\TDA\big)$ is dense in $L^2(\cA)$ by Lemma
\ref{lem.3}, we can construct a sequence $\{f_n\}_{n \ge 1}$ in
$H^{-\frac{1}{2}}(\GD)$ such that
\begin{equation}
\label{eq:PF1}
\|\TDA f_n - G(\cdot,\vz)\|_{L^2(\cA)} < \frac{1}{n}, \qquad n \ge 1.
\end{equation}
Lemma \ref{lem.3} also states that $\mbox{range}\big({\TAD}\big)$ is
dense in $H^{-\frac{1}{2}}(\GD)$, so we can construct a sequence $\{g_n\}_{n
  \ge 1}$ in $L^2(\cA)$ satisfying
\begin{equation}
\|\TAD g_n - f_n\|_{H^{-\frac{1}{2}}(\GD)} < \frac{1}{n}, \qquad n \ge 1.
\label{eq:PF2}
\end{equation}
These results, the triangle inequality and Lemma \ref{lem.1} give 
\begin{align}
\|N g_n - G(\cdot,\vz)\|_{L^2(\cA)} &= \|\TDA \TAD g_n - G(\cdot,
\vz)\|_{L^2(\cA)} \nonumber \\ &\le \big\|\TDA \big(\TAD g_n - f_n)
\big\|_{L^2(\cA)} + \|\TDA f_n - G(\cdot,\vz)\|_{L^2(\cA)} \nonumber \\ &<
\frac{\|\TDA\| + 1}{n}. \label{eq:PF6}
\end{align}
By the Archimedian property of real numbers, $\forall \, \ep >0$,
there exists a natural number $N$ such that $\big(\|\TDA\| + 1\big)/n
< \ep$, for all $n > N$, so we have shown that \eqref{eq:1.15} holds.

It remains to prove that the sequence $\Big\{ \|\TAD
g_n\|_{H^{-\frac{1}{2}}(\GD)}\Big\}_{n \ge 1}$ cannot be bounded. We argue
by contradiction: Suppose that this sequence were bounded. Then, we
obtain from \eqref{eq:PF2} that
$\{\|f_n\|_{H^{-\frac{1}{2}}(\GD)}\}_{n \ge 1}$ is a bounded sequence,
so there exists a subsequence $\{f_{n_m}\}_{m \ge 1}$ that converges
weakly to some $f \in H^{-\frac{1}{2}}(\GD)$. By \eqref{eq:PF2} this
means
\begin{equation}
\TAD g_{n_m} \to f, \qquad \mbox{weakly in }H^{-\frac{1}{2}}(\GD),
\label{eq:PF3}
\end{equation}
and since $\TDA$ is compact by Lemma \ref{lem.3}, we have
\begin{equation}
N g_{n_m} = \TDA \TAD g_{n_m} \to \TDA f, \qquad \mbox{strongly in }
L^2(\cA).
\label{eq:PF4}
\end{equation}
But \eqref{eq:PF1} implies that $\TDA f = G(\cdot, \vz)\big|_{\cA}$, which
contradicts \eqref{eq:PF0}. This proves that the sequence $\Big\{ \|\TAD
g_n\|_{H^{-\frac{1}{2}}(\GD)}\Big\}_{n \ge 1}$ cannot be bounded, as stated in the
theorem. \endproof

\begin{remark}\label{rem1}
The statement of Theorem \ref{thm.1}, which is based on the validity
of Lemmas \ref{lem.1}--\ref{lem.3}, holds for any  {wave number} $k \in \RR$ with
the exception of a discrete set of isolated values. These exceptional
points correspond to either $-k^2$ being a  {Neuman} eigenvalue of the Laplacian in
$\cD$ or to values of $k^2$ at which the forward problem
\eqref{eq:1.7}--\eqref{eq:1.8} is not uniquely solvable.  More details
are in appendix \ref{ap:A}.
\end{remark}

\subsubsection{The imaging algorithm}
\label{sect:ImagAlg}
Suppose that the imaging region is the sector $(x_{_I},0) \times \cX$ of $\cW_o$, 
with $x_{_I} > \xA$ satisfying
\begin{equation}
x_{_I} - \xA > \frac{1}{|\beta_{J+1}|},
\label{eq:IMAs}
\end{equation}
so that we can neglect all the evanescent modes. 
Using the mode decomposition of the scattered wave, we can rewrite
\eqref{eq:1.15} as a linear least squares problem for a $(J+1) \times
(J+1)$ linear system of equations.  Indeed, by linear superposition,
we can  decompose the scattered field as
\begin{equation}
u(\vx_r,\vx_s) = \sum_{j=0}^\infty u_j(\vx_r) \psi_j(\bxp_s),
\label{eq:IMP1}
\end{equation}
where $u_j(\vx)$ solves \eqref{eq:1.7}--\eqref{eq:1.8}, with
$G(\vx,\vx_s)$ replaced in \eqref{eq:1.7b} by
\[
G_j(\vx) = \int_{\cX} d \bxp_s \, G(\vx,\vx_s) \psi_j(\bxp_s).
\]  
Furthermore, we can represent the array response  \eqref{eq:1.9}
by the $(J+1) \times (J+1)$ matrix ${\bf U} = \big(U_{j,j'}\big)_{0
  \le j \le J}$ with entries
\begin{equation}
U_{j,j'} = \int_{\cX} d \bxp_r \int_{\cX} d \bxp_s \, u(\vx_s,\vx_r)
\psi_j(\bxp_r) \psi_{j'}(\bxp_s),
\label{eq:IMP2}
\end{equation}
where we recall the assumption \eqref{eq:FA}.

Neglecting the evanescent modes, we obtain from the definition
\eqref{eq:1.10} of the near field operator {that} 
\begin{equation}
u(\vx_r,\vx_s) \approx \sum_{j,j' = 0}^J  U_{j,j'}
\psi_j(\bxp_r) \psi_{j'}(\bxp_s),
\label{eq:IMP2p}
\end{equation}
{and} 
\begin{equation}
N g(\vx_r) \approx \sum_{j=0}^J\psi_{j}(\bxp_r) \sum_{j'=0}^J U_{j,j'}
g_{j'} = \sum_{j=0}^J\psi_{j}(\bxp_r) \big({\bf U}{\bf g}\big)_j, \qquad \forall \, \vx_r \in \cA,
\label{eq:IMP3}
\end{equation}
where ${\bf g} = (g_o, \ldots,
g_J)^T$ is the $J+1$ column vector with components
\begin{equation}
g_j = \int_{\cX} d \bxp_s \, \psi_j(\bxp_s) g(\vx_s).
\label{eq:IMP4}
\end{equation}
Moreover, using the assumption \eqref{eq:IMAs}, 
\begin{equation}
G(\vx_r,\vz) \approx \sum_{j=0}^J b_{j,\vz} \, \psi_j(\bxp_r), 
\label{eq:IMP4p}
\end{equation}
with 
\begin{equation}
b_{j,\vz} = \int_{\cX} d \bxp \, \psi_j(\bxp) G(\vx,\vz), \qquad \vx =
(\xA,\bxp) \in \cA.
\label{eq:IMP5}
\end{equation}
Letting ${\bf b}_{\vz}$ be the $J+1$ column vector with components 
\eqref{eq:IMP5},  we obtain that
\begin{equation}
N g(\vx_r) - G(\vx_r,\vz) \approx \sum_{j=0}^J \psi_j(\bxp_r) \Big( {\bf
  U}{\bf g} - {\bf b}_{\vz} \big)_j, \qquad \forall \, \vx_r \in \cA.
\label{eq:IMP6}
\end{equation}
The eigenfunction are orthonormal, so we can write
\begin{equation}
\|N g - G(\cdot, \vz)\|_{L^2(\cA)} \approx \|{\bf U}{\bf g} - {\bf
  b}_{\vz}\|_2,
\label{eq:IMP7}
\end{equation}
where $\| \cdot \|_2$ is the Euclidian norm.
The summary of the linear sampling algorithm for estimating $\GD$ is
as follows:

\vspace{0.05in}
\begin{algorithm}
\label{alg.1}

\noindent \textbf{Input:} The $(J+1)\times(J+1)$ matrix ${\bf U}$ and the imaging mesh.

\noindent \textbf{Processing steps:}
\vspace{0.05in}
\begin{enumerate}
\item For a user defined small $\ep>0$, and for all $\vz $ on the imaging mesh, solve the normal equations
\begin{equation}
\big({\bf U}^\star {\bf U} + \alpha^\ep {\bf I} \big) {\bf g}_{\vz} =
    {\bf U}^\star {\bf b}_{\vz},
\label{eq:IMP8}
\end{equation}
where ${\bf U}^\star$ is the Hermitian adjoint of ${\bf U}$, ${\bf I}$
is the $(J+1) \times(J+1)$ identity matrix and $\alpha^\ep$ is a
positive Tikhonov regularization parameter chosen according to the
Morozov principle, so that
\[
\|{\bf U}{\bf g}_{\vz} - {\bf b}_{\vz}\|_2 = \ep \|{\bf g}_{\vz}\|_2.
\]
\item Calculate the indicator function
\begin{equation}
{\cal J}(\vz) = \frac{1}{\|{\bf g}_{\vz}\|_2}.
\label{eq:IMP9}
\end{equation}
\end{enumerate}

\noindent \textbf{Output:} The estimate of the support of $\cD$ is
determined by the set of points $\vz$ where ${\cal J}(\vz)$ exceeds a
user defined threshold. The estimated wall deformation $\GD$ is the
part of the boundary of $\cD$ contained in $\cW_o$.
\end{algorithm}

\subsection{Imaging with a partial aperture array}
\label{sect:Prob1.3}
If the  array does not cover the entire cross-section of the waveguide,
\begin{equation}
\cA = \{\xA\} \times \cX_{\cA}, \qquad \cX_{\cA} \subset \cX,
\label{eq:PA1}
\end{equation}
we can calculate the analogue of \eqref{eq:IMP2}, the $(J+1) \times
(J+1)$ matrix ${\bf U}^{\cA} = ({U}^{\cA}_{j,j'})_{0
  \le j,j'\le J}$ with entries
\begin{equation}
U^{\cA}_{j,j'} = \int_{\cX_{\cA}} d \bxp_r \int_{\cX_{\cA}} d \bxp_s \, u(\vx_s,\vx_r)
\psi_j(\bxp_r) \psi_{j'}(\bxp_s), \qquad j,j' = 0, \ldots, J.
\label{eq:PA2}
\end{equation}
This is related to ${\bf U}$ by
\begin{equation}
{\bf U}^{\cA} \approx {\bf M} {\bf U} {\bf M},
\label{eq:PA3}
\end{equation}
where we used the approximation \eqref{eq:IMP2p} and introduced the
symmetric, positive semidefinite Gram matrix ${\bf M} = (M_{j,j'})_{0
  \le j,j'\le J}$ with entries
\begin{equation}
M_{j,j'} = \int_{\cX_{\cA}}d \bxp \, \psi_j(\bxp) \psi_{j'}(\bxp),
\qquad j,j' = 0, \ldots, J.
\label{eq:PA4}
\end{equation}
While ${\bf M}$ equals the identity when the array has full aperture,
at partial aperture it is poorly conditioned. Thus, we cannot
calculate ${\bf U}$ from \eqref{eq:PA3} by
inverting the Gramian ${\bf M}$.
If we let
\begin{equation}
{\bf M} = {\bf V} \mbox{diag}(\sigma_0, \ldots, \sigma_J) {\bf V}^T,
\label{eq:PA5}
\end{equation}
be the eigenvalue decomposition of ${\bf M}$, with ${\bf V} =
(\bv_j)_{0\le j \le J}$ the orthogonal matrix of eigenvectors $\bv_j$,
and with the eigenvalues in decreasing order $\sigma_o \ge \sigma_1 \ge \ldots
\ge \sigma_J \ge 0, $ then we expect that 
\begin{equation}
0 \le  \sigma_j \ll 1, \qquad J_M < j \le  J,
\label{eq:PA7}
\end{equation}
for some $J_M < J$. Then, we 
approximate ${\bf U}$ by
\begin{equation}
\widetilde{{\bf U}} = {\bf M}^\dagger {\bf U}^{\cA} 
{\bf M}^\dagger \approx  {\bf M}^\dagger {\bf M} {\bf U} {\bf M} {\bf M}^\dagger,
\label{eq:PA8}
\end{equation}
with
\begin{equation}
{\bf M}^\dagger = {\bf V} \mbox{diag}(\sigma_0^{-1}, \ldots,
\sigma_{J_M}^{-1}, 0, \ldots, 0) {\bf V}^T.
\label{eq:PA9}
\end{equation}
Note that ${\bf M}^\dagger {\bf M}$ is the orthogonal projection on
$\mbox{span}(\bv_o, \ldots, \bv_{J_M})$.

The imaging algorithm is almost the same as Algorithm \ref{alg.1},
except that the input matrix is replaced by $\widetilde{\bf U}$, which
we can compute, and ${\bf b}_{\vz}$ is replaced by
\begin{equation}
\widetilde{\bf b}_{\vz} = {\bf M}^\dagger {\bf M} {\bf b}_{\vz}.
\label{eq:PA8p}
\end{equation}
To give a more concrete explanation of the effect of the aperture, let
us use definition \eqref{eq:IMP2} and equation \eqref{eq:PA8} to
relate $\widetilde{\bf U}$ to the full aperture response
\begin{equation}
\widetilde{\bf U} = \sum_{j,j'=0}^{J_M} \bv_j \bv_{j'}^T
\int_{\cX} d \bxp_r \, p_j(\bxp_s)\int_{\cX} d \bxp_s\, 
p_{j'}(\bxp_s)u(\vx_r,\vx_s),
\label{eq:PA11}
\end{equation}
where now $\vx_r,\vx_s \in \{\xA\} \times \cX$ and 
\begin{equation}
p_j(\bxp) = \sum_{l=0}^J v_{l,j} \psi_l(\bxp).
\label{eq:PA12}
\end{equation}
The vector \eqref{eq:PA8p} is 
\begin{equation}
\widetilde{\bf b}_{\vz} = \sum_{j=0}^{J_M} \bv_j \int_{\cX} d \bxp_r
\, p_j(\bxp_r) G\big(\vx_r,\vz\big),
\label{eq:PA11p}
\end{equation}
and if we use ${\bf g}$ defined in \eqref{eq:IMP4}, we obtain 
\begin{align}
\| \widetilde{\bf U} {\bf g} -\widetilde{\bf b}_{\vz}\|_2^2 =
\sum_{j=0}^{J_M}& \left| \int_{\cX} d \bxp_r p_j(\bxp_r) \left[
  \int_{\cX} d \bxp_s\, u(\vx_r,\vx_s) \widetilde{g}\big(\vx_s) -
  G\big(\vx_r,\vz\big)\right] \right|^2,
\label{eq:PA20}
\end{align} 
with 
\begin{equation}
\widetilde g(\vx_s)= \sum_{j = 0}^{J_M} p_j(\byp) \bv_j^T {\bf g} =
\sum_{j = 0}^{J_M} p_j(\bxp_s) \int_{\cX} d \bxp \, p_j(\bxp)
g\big((\xA,\bxp)\big).
\label{eq:PA21}
\end{equation}
We can also define the analogue of \eqref{eq:IMP2p}
\begin{align}
\widetilde u(\vx,\vy) &= \sum_{l,l'=0}^{J} \widetilde U_{l,l'}
\psi_l(\bxp) \psi_{l'}(\byp) \nonumber \\ &= \sum_{j,j'=0}^{J_M}
p_j(\bxp) p_{j'}(\byp) \int_{\cX} d \bxp_r \, p_j(\bxp)\int_{\cX} d
\bxp_{s}\, p_{j'}(\bxp_s)u(\vx_r,\vx_s),
\label{eq:PA22}
\end{align}
for all $\vx,\vy \in \{\xA\} \times \cX.$

Note that $\{p_j(\bxp)\}_{0 \le j \le J}$ is an orthogonal set in $\mbox{span} 
\{ \psi_j(\bxp), ~ 0 \le j \le J\}$, satisfying 
\begin{equation}
\label{eq:PA13}
\int_{\cX} d \bxp \, p_j(\bxp) p_{j'}(\bxp) = \delta_{j,j'}, \qquad
\int_{\cX_{\cA}} d \bxp \, p_j(\bxp) p_{j'}(\bxp) = \sigma_j
\delta_{j,j'}.
\end{equation}
Therefore, $\widetilde u$, $\widetilde{\bf b}_{\vz}$ and
$\widetilde{g}$ are projections of their continuum aperture
counterparts on the subspace $\mbox{span} \{ p_j(\bxp), ~ 0 \le j \le
J_M\}$.  The second relation in \eqref{eq:PA13} shows that $\sigma_j
\in [0,1]$ and we must have
\begin{equation}
p_j(\bxp) \approx 0 ~~\mbox{in}~ \cX_{\cA}, ~~\mbox{for}~~j = J_M+1,
\ldots, J,
\label{eq:PA14}
\end{equation}
and 
\begin{equation}
p_j(\bxp) \approx 0 ~~\mbox{in}~ \cX \setminus \cX_{\cA},
~~\mbox{for}~~ \sigma_j \approx 1.
\label{eq:PA15}
\end{equation}
We verify in the next section, for a two dimensional waveguide, that
$\sigma_j\approx 1$ for $0 < j < J_M$, where $J_M = \lfloor J
{|\cX_{\cA}|}/{|\cX|}\rfloor$ and $|\cX_{\cA}|$, $|\cX|$ are the
lengths of the aperture and cross-section of the waveguide.  Thus, the
projection limits the support of the functions to the array aperture
$\cA$.

\subsubsection{Illustration in a two dimensional waveguide}
In two dimensions, the cross-section of the waveguide is the interval
$\cX = (0,|\cX|)$ of length $|\cX|$. Suppose that the array aperture
is \[\cA= \{\xA\} \times \cX_{\cA}, \qquad \cX_{\cA} =
(0,|\cX_{\cA}|), \qquad |\cX_{\cA}| < |\cX|.
\] 
Then, using the
eigenfunctions \eqref{eq:1.3} of the Laplacian
\begin{equation}
\psi_0(\bxp) = \frac{1}{\sqrt{|\cX|}}, \qquad \psi_j(\bxp) =
\sqrt{\frac{2}{|\cX|}} \cos \left( \frac{ \pi j \bxp}{|\cX|} \right),
\qquad j \ge 1,
\label{eq:PA30}
\end{equation}
we obtain that the Gram matrix ${\bf M}$ is 
\begin{equation}
M_{j,j'} = \left\{ \begin{array}{ll} \frac{|\cX_{\cA}|}{|\cX|}, \qquad
  &j=j'=0, \\ \frac{|\cX_{\cA}|}{|\cX|} \sqrt{2} \,\mbox{sinc}
  \left(\pi j' \frac{|\cX_{\cA}|}{|\cX|}\right), & j=0, ~ 1 \le j' \le
  J, \\ \frac{|\cX_{\cA}|}{|\cX|} \sqrt{2} \,\mbox{sinc} \left(\pi j
  \frac{|\cX_{\cA}|}{|\cX|}\right), & j'=0, ~ 1 \le j \le J,
  \\ \frac{|\cX_{\cA}|}{|\cX|} \left[ \mbox{sinc}
    \left(\pi(j-j')\frac{|\cX_{\cA}|}{|\cX|}\right) + \mbox{sinc}
    \left(\pi(j+j')\frac{|\cX_{\cA}|}{|\cX|}\right) \right], & 1 \le
  j,j', \le J.
\end{array} \right.
\end{equation}

The eigenvalues of ${\bf M}$ are related to the eigenvalues of the
$(2J+1) \times (2J+1)$ prolate matrix \cite{varah,slepian}, which is
symmetric and Toeplitz
\begin{equation}
{\bf T} = \begin{pmatrix} 
t_0 & t_1 & t_2  & \ldots & t_{2J} \\
t_1 & t_0 & t_1 &\ldots & t_{2J-1} \\
t_2 & t_1 & \ddots & \ddots & \vdots \\
\vdots & & \ddots & \ddots & t_1 \\
t_{2J} & \ldots& \ldots & t_1& t_0.
\end{pmatrix}, \qquad 
t_j = \frac{|\cX_{\cA}|}{|\cX|} \mbox{sinc} \left(\pi j
\frac{|\cX_{\cA}|}{|\cX|}\right).
\end{equation}
To make the connection to ${\bf M}$, we rewrite ${\bf T}$ as the
matrix
\begin{equation}
{T}_{j,j'} = \frac{|\cX_{\cA}|}{|\cX|} \mbox{sinc} \left(\pi (j-j')
\frac{|\cX_{\cA}|}{|\cX|}\right), \qquad -J \le j,j' \le J,
\end{equation}
using that $t_{j-j'} = T_{j,j'}$, for $-J \le j' \le j \le J$. This
matrix has $J$ odd eigenvectors $\{\boldsymbol{\tau}_j^{o}\}_{1 \le j
  \le J}$ for eigenvalues $\{\sigma_j^o\}_{1 \le j \le J}$ and $J+1$
even eigenvectors $\{\boldsymbol{\tau}_j^{e}\}_{0 \le j \le J}$ for
eigenvalues $\{\sigma_j^e\}_{0\le j \le J}$. Odd and even means that
the components $\tau_{l,j}^o$ and $\tau_{l,j}^e$ of the eigenvectors
satisfy
\[
\tau_{-l,j}^o = - \tau_{l,j}^o, \qquad \tau_{-l,j}^e = \tau_{l,j}^e,
\qquad l = 1, \ldots, J.
\]

We are interested in the even spectrum of ${\bf T}$,  which determines the eigenvalues  
$\sigma_j = \sigma_j^e$ of ${\bf M}$, with the eigenvectors given by 
\begin{equation}
\bv_j = (v_{0,j}, \ldots, v_{J,j})^T, \qquad 
v_{l,j} = \left\{ \begin{array}{ll}
\sqrt{2} \, \tau_{0,j}^e, \qquad &l = 0 \\
\tau_{l,j}^e, & 1 \le l \le J.
\end{array} \right.
\end{equation}
Then, we conclude from the known properties \cite{slepian}  of the spectrum of ${\bf T}$ 
that $ \sigma_j \approx 1$ for
$0 \le j < J_M = \left \lfloor J \frac{|\cX_{\cA}|}{|\cX|} \right
\rfloor,$ and that $ \sigma_j \approx 0$ for $ j > J_M.$ Moreover, the
orthogonal functions $p_j(\bxp)$ defined in \eqref{eq:PA12} are
trigonometric polynomials supported in $\cX_{\cA}$ for $0 \le j < J_M$
and in $\cX \setminus \cX_{\cA}$ for $j > J_M$, as stated in the
previous section. At the threshold index $j = J_M$, the polynomial
$p_{J_M}(\bxp)$ is sharply peaked at the end of the interval $\cX_{\cA}$
\cite{slepian}.

\section{Imaging inside the waveguide with wall deformations}
\label{sect:Problem2}
The analysis of the linear sampling method for estimating both the
support $\Omega$ of scatterers in the waveguide and the wall
deformation $\Gamma$ is very similar to that in the previous section,
so we do not include it here and state directly the results.

The near field operator is defined as in \eqref{eq:1.10}, using the
scattered wave $u(\vx_r,\vx_s)$ at the array, and its factorization is
similar to \eqref{eq:1.11}
\begin{equation}
N = \TDOA \TADO,
\label{eq:2.5}
\end{equation}
where the operators $\TDOA$ and $\TADO$ are the analogues of $\TDA$
and $\TAD$ defined in Lemma \ref{lem.1}.

In the case of an impenetrable scatterer, the field
$u(\vx,\vx_s)$ satisfies
\begin{align}
\big(\Delta_{\vx} + k^2 \big) u(\vx,\vx_s) &= 0, \qquad \vx \in \cW
\setminus { \overline{\Omega}},
\label{eq:2.1}\\
 \frac{\partial u(\vx,\vx_s)}{\partial \vnu} &= 0, \quad \vx \in
  {  \partial W_o \setminus \overline{\Gamma}_o}, \\ \frac{\partial
   u(\vx,\vx_s)}{\partial \vnu} &= - \frac{\partial
   G(\vx,\vx_s)}{\partial \vnu}, \qquad \vx \in \Gamma, 
\label{eq:2.2}\\
\mathcal{B} u(\vx,\vx_s) &= - \mathcal{B} G(\vx,\vx_s), \qquad \vx \in
\partial \Omega,
\end{align}
and the radiation condition in Definition \ref{def.1}, where
$\mathcal{B} u = u$ if the scatterer is sound soft and $\mathcal{B} u
= \partial_{\vnu} u$ if it is sound hard ({or more generally $\mathcal{B} u$ maybe be a combination of Robin type}).  

For a penetrable scatterer, modeled by the square $n^2(\vx)$ of the
index of refraction, with positive real part $\Re(n^2) > 0$ and
non-negative imaginary part $\Im(n^2) \ge 0$, and with support of
$n^2(\vx) - 1$ in $\Omega$, the scattered field satisfies
\begin{align}
\big(\Delta_{\vx} + k^2 n^2(\vx)\big) u(\vx,\vx_s) &= -
k^2\big(n^2(\vx)-1) G(\vx,\vx_s), \qquad \vx \in \cW,
\label{eq:2.3}\\
 \frac{\partial u(\vx,\vx_s)}{\partial \vnu} &= 0, \quad \vx \in
 {  \partial W_o \setminus \overline{\Gamma}_o}, \\ \frac{\partial
   u(\vx,\vx_s)}{\partial \vnu} &= - \frac{\partial
   G(\vx,\vx_s)}{\partial \vnu}, \qquad \vx \in \Gamma, 
\label{eq:2.4}
\end{align}
and the radiation condition  in Definition \ref{def.1}.

The operators $\TDOA$ and $\TADO$ are defined as in Lemma \ref{lem.1},
with $\GD$ replaced by $\GD \cup \partial \Omega$, when the scatterer
is sound hard.  

For a sound soft scatterer we define $ \TADO:L^2(\cA) \to
H^{-\frac{1}{2}}(\GD) \times H^{\frac{1}{2}}(\partial \Omega) $ by
\begin{equation}
\TADO g(\vz,\vz') = \left(\partial_{\vnuz} \hspace{-0.03in} \int_{\cA} d S_{\vx_s}\,
G(\vz,\vx_s) g(\vx_s), \int_{\cA} d S_{\vx_s}\,
G(\vz',\vx_s) g(\vx_s)\right), \label{eq:2.10}
\end{equation}
for arbitrary points $\vz \in \GD$ and $\vz' \in \partial \Omega$ and
for arbitrary $g \in L^2(\cA)$. The operator
$\TDOA:H^{-\frac{1}{2}}(\GD) \times H^{\frac{1}{2}}(\partial \Omega)
\to L^2(\cA)$ takes arbitrary functions $f_{\GD} \in
H^{-\frac{1}{2}}(\GD)$ and $f_{\partial \Omega} \in
H^{\frac{1}{2}}(\partial \Omega)$ and returns the trace
$\TDOA(f_{\GD},f_{\partial \Omega}) = w\big|_{\cA}$ of the solution of
\begin{align}
\big(\Delta_{\vx} + k^2 \big) w(\vx) &= 0, \qquad \vx \in \cW
\setminus {\overline{\Omega}}, \label{eq:2.11} \\ \frac{\partial w(\vx)}{\partial
  \vnu} &= 0, \quad \vx \in  {  \partial W_o \setminus \overline{\Gamma}_o},
\\ \frac{\partial w(\vx)}{\partial \vnu} &= - f_{\GD}(\vx), \qquad \vx
\in \Gamma, \label{eq:2.12} \\ w(\vx) &= - f_{\partial \Omega}(\vx),
\qquad \vx \in \partial \Omega,
\end{align}
satisfying the radiation condition as in Definition \ref{def.1}.

For a penetrable scatterer, the operator $ \TADO:L^2(\cA) \to
H^{-\frac{1}{2}}(\GD) \times H^{1}(\Omega) $ is
\begin{equation}
\TADO g(\vz,\vz') = \left(\partial_{\vnuz} \hspace{-0.03in} \int_{\cA} d S_{\vx_s}\,
G(\vz,\vx_s) g(\vx_s), \int_{\cA} d S_{\vx_s}\,
G(\vz',\vx_s) g(\vx_s)\right), \label{eq:2.14}
\end{equation}
for arbitrary points $\vz \in \GD$, $\vz' \in \Omega$ and functions $g
\in L^2(\cA)$. Moreover, the operator $\TDOA:H^{-\frac{1}{2}}(\GD)
\times H^{1}(\Omega) \to L^2(\cA)$ is defined by the trace
$\TDOA(f_{\GD},f_{\Omega}) = w\big|_{\cA}$ of the solution of the
boundary value problem
\begin{align}
\big(\Delta_{\vx} + k^2n^2(\vx) \big) w(\vx) &= -k^2
\big(n^2(\vx)-1\big) f_{\Omega}(\vx), \qquad \vx \in
\cW, \label{eq:2.15} \\ \frac{\partial w(\vx)}{\partial \vnu} &= 0,
\quad \vx \in  {  \partial W_o \setminus \overline{\Gamma}_o}, \\ \frac{\partial
  w(\vx)}{\partial \vnu} &= - f_{\GD}(\vx), \qquad \vx \in
\Gamma, \label{eq:2.16}
\end{align}
satisfying a radiation condition as in Definition \ref{def.1}, for
arbitrary $f_{\GD} \in H^{-\frac{1}{2}}(\GD)$ and $f_{ \Omega} \in
H^{1}(\Omega)$.

The analogue of Theorem \ref{thm.1} is: 

\vspace{0.05in}
\begin{theorem}
\label{thm.2}
Let $ \vz $ be a search point in $\cW_o$,
between the array and the end wall. For any $\ep > 0$ {let $g_{\vz}^\ep \in L^2(\cA)$ satisfy}
\begin{equation}
\|Ng_{\vz}^\ep - G(\cdot,\vz)\|_{L^2(\cA)} < \ep.
\label{eq:2.20}
\end{equation}
Let $H$ denote $H^{-\frac{1}{2}}(\GD\cup \partial \Omega)$ in the case
of a sound hard scatterer, or $H^{-\frac{1}{2}}(\GD)\cup
H^{\frac{1}{2}}(\partial \Omega)$ for a sound soft scatterer, or
$H^{-\frac{1}{2}}(\GD)\cup H^{1}(\Omega)$ for a penetrable scatterer.
There are two possibilities:

\vspace{0.03in}
\begin{enumerate}
\item If $\vz \in \cD \cup \Omega$, {there exists a} $g_{\vz}^\ep$ satisfying
  \eqref{eq:2.20}, {such that}    $\|\TADO g_{\vz}^\ep\|_{H}$
  remains bounded as $\ep \to 0$.
\item If $\vz \notin \cD\cup \Omega$, for any $g_{\vz}^\ep$ satisfying
  \eqref{eq:2.20}, $\displaystyle \lim_{\ep \to 0} \|\TADO
  g_{\vz}^\ep\|_{H} = \infty$.
\end{enumerate}
\end{theorem}
\begin{remark}
The statement of Theorem \ref{thm.2} also holds for any wave number $k \in \RR$ with
the exception of a discrete set of isolated values. In this case, in addition to the exceptional
wave numbers in Remark \ref{rem1}, one has to exclude the values of $k$ for which $-k^2$ is an eigenvalue of the Laplacian in $\Omega$ with the respective boundary condition in the case of impenetrable scatterer or a transmission eigenvalue in  $\Omega$ in the case of penetrable scatterer (for the latter see \cite{CCH}). 
\end{remark}

\vspace{0.05in} As in the previous section, the imaging is based on
the indicator function ${1}/{\|{\bf g}_{\vz}\|_2}$, which is expected
to be very small for points $\vz \notin \cD \cup \Omega$. Algorithm
\ref{alg.1} remains unchanged, which is useful because in practice it
is not known if the waveguide is empty or not. The case of a partial
aperture array is handled the same way as in section
\ref{sect:Prob1.3}.

\section{Numerical results}
\label{sect:Numerics}
We assess the performance of the linear sampling algorithm using
numerical simulations in a two dimensional waveguide.  All the
coordinates are scaled by the width $|\cX|$ of the waveguide, and we
vary the wavelength to get a smaller or larger number of propagating modes
\[
J +1 = \left\lfloor k
\frac{|\cX|}{\pi}\right \rfloor+ 1.\] 
The array data $u(\vx_r,\vx_s)$ are obtained by solving the wave
equation in the sector $(-8|\cX|,0) \times (0,|\cX|)$ of the waveguide, using
the high performance multi-physics finite element software
Netgen/NGSolve \cite{schoberl1997netgen} and a perfectly matched layer
at the left end of the domain. The separation between the sensors is of the order 
of the wavelength, more precisely: $\frac{\cX}{25}$ in the case of $10$ and $20$ propagating modes, and  $\frac{\cX}{55}$ in the case of $50$ propagating modes.  The matrices ${\bf U}$ and ${\bf U}^{\cA}$ are 
calculated as in \eqref{eq:IMP2} and \eqref{eq:PA2}, by approximating the 
integrals with  Simpson's quadrature rule. The data are contaminated with 2\% multiplicative noise, meaning that   the $ij$-th entry of the contaminated matrix is  $ {U}_{ij} (1+ 0.02 \delta \ )$ where $\delta$ is a uniformly distributed random number between $0$ and $1$.

The imaging region is $(-4|\cX|,0)\times (0,|\cX|)$ and the array is at range { $\xA = - 5 |\cX|$}. The
images are obtained with Algorithm \ref{alg.1} in the case of a full
aperture $\cA = \{\xA\} \times (0,|\cX|)$ or its modification
explained in section \ref{sect:Prob1.3} in the case of partial
aperture $\cA = \{\xA\} \times (0,|\cX_{\cA}|)$, with $|\cX_{\cA}| <
|\cX|$.  For better visualization we display the logarithm of the
indicator function \eqref{eq:IMP9}.

\begin{figure}[t]
\begin{center}
\includegraphics[scale=0.25]{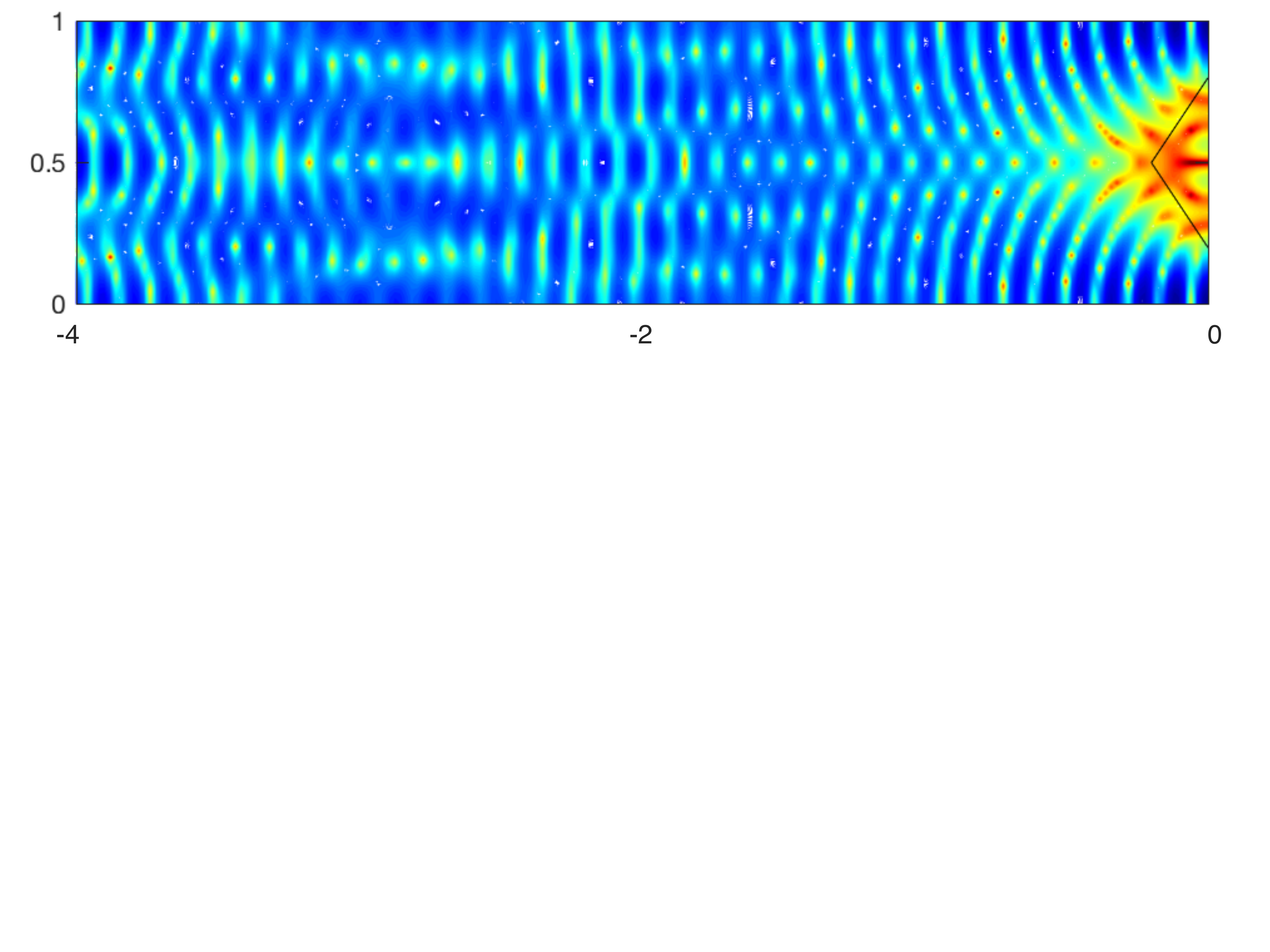}
\vspace{-4cm}
\end{center}
\begin{center}
\includegraphics[scale=0.25]{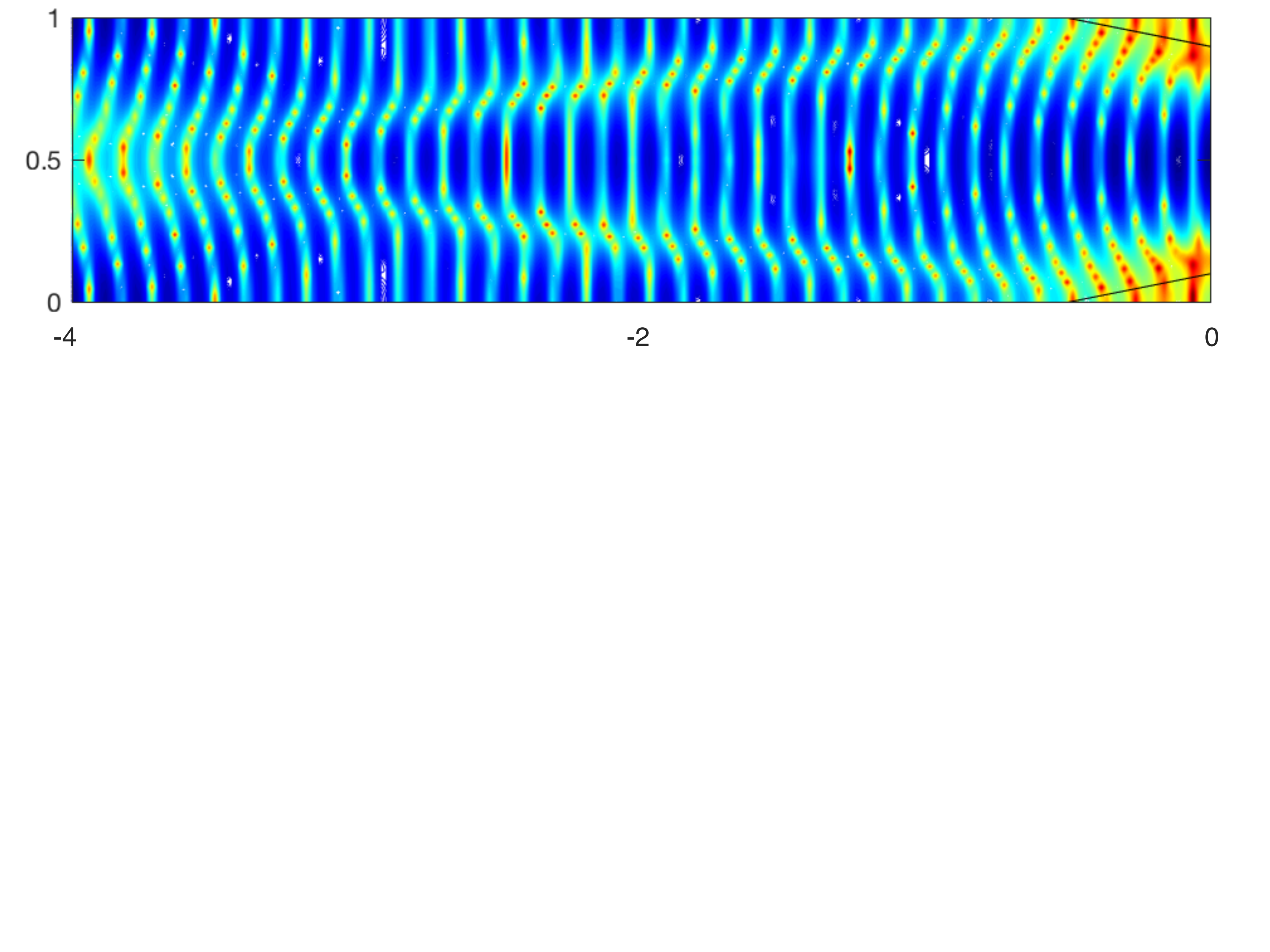}
\vspace{-4cm}
\caption{Reconstruction of wall deformations shown with a solid black
  line.  The abscissa is range and the ordinate is cross-range scaled
  by $|\cX|$. Full aperture data and { $J +1= 10$}. }
\label{side1_10modes}
\end{center}
\end{figure}

\begin{figure}[H]
\begin{center}
\includegraphics[scale=0.25]{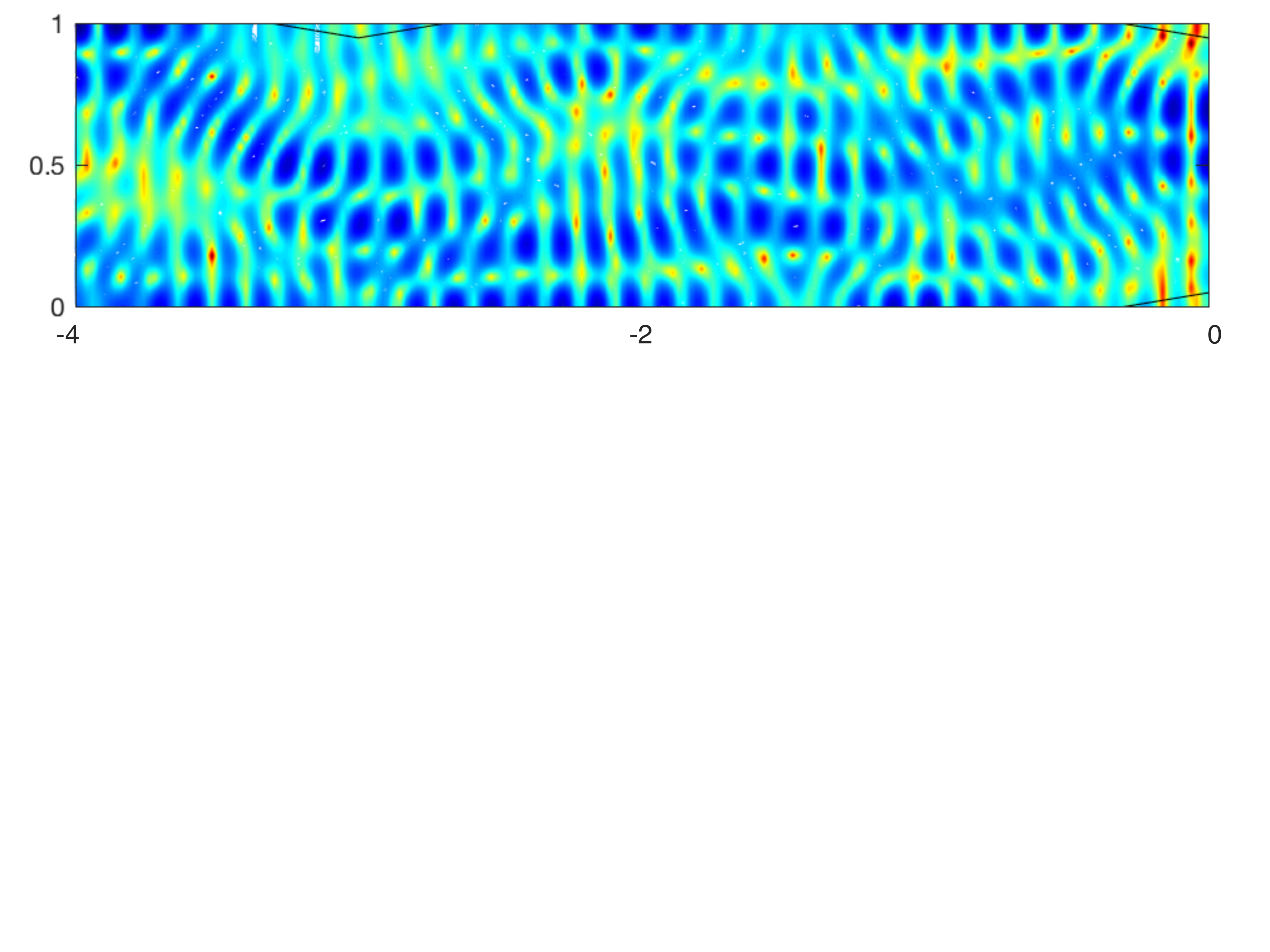}
\vspace{-4cm}
\end{center}
\begin{center}
\includegraphics[scale=0.25]{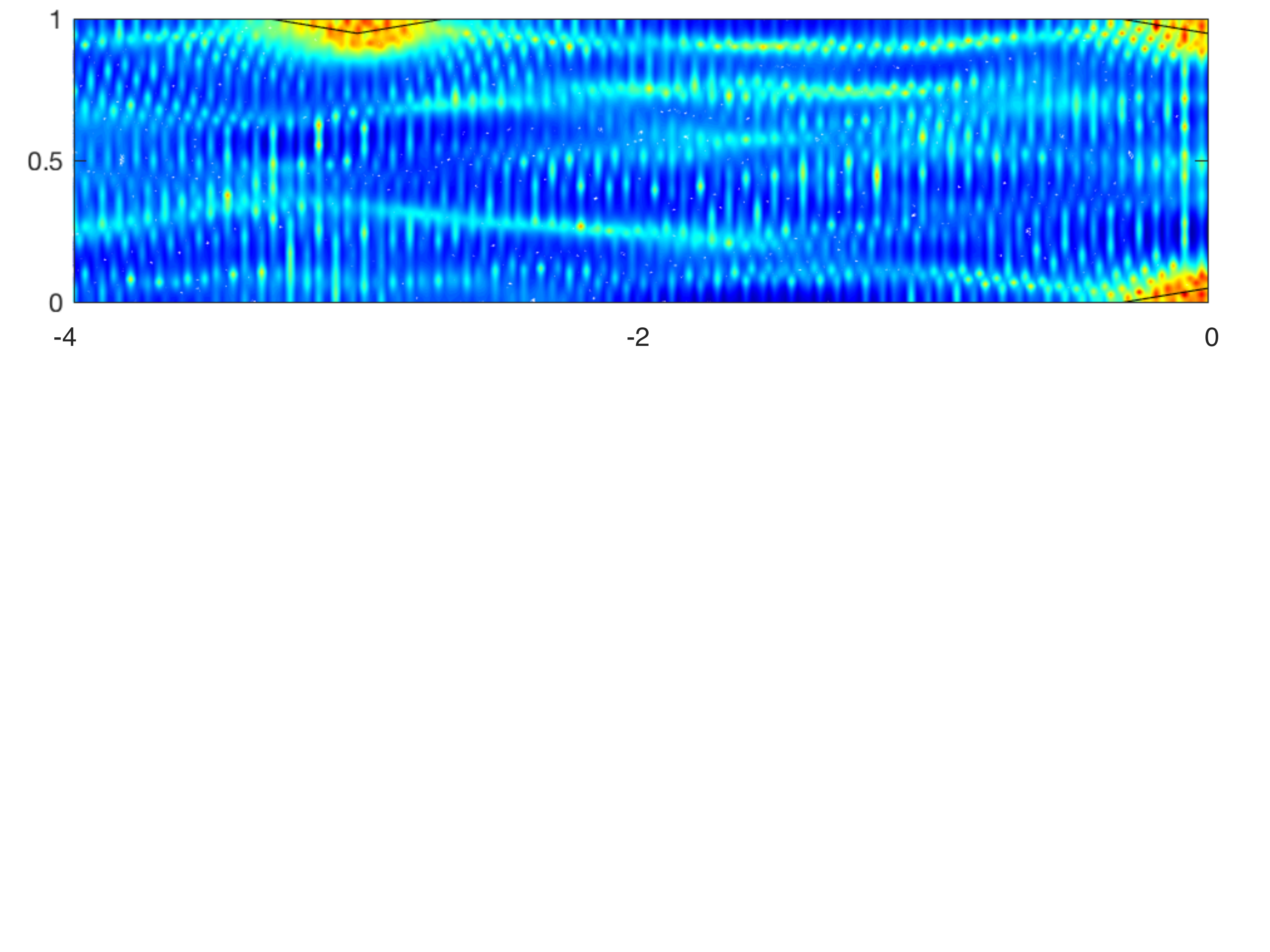}
\vspace{-4cm}
\end{center}
\begin{center}
\includegraphics[scale=0.25]{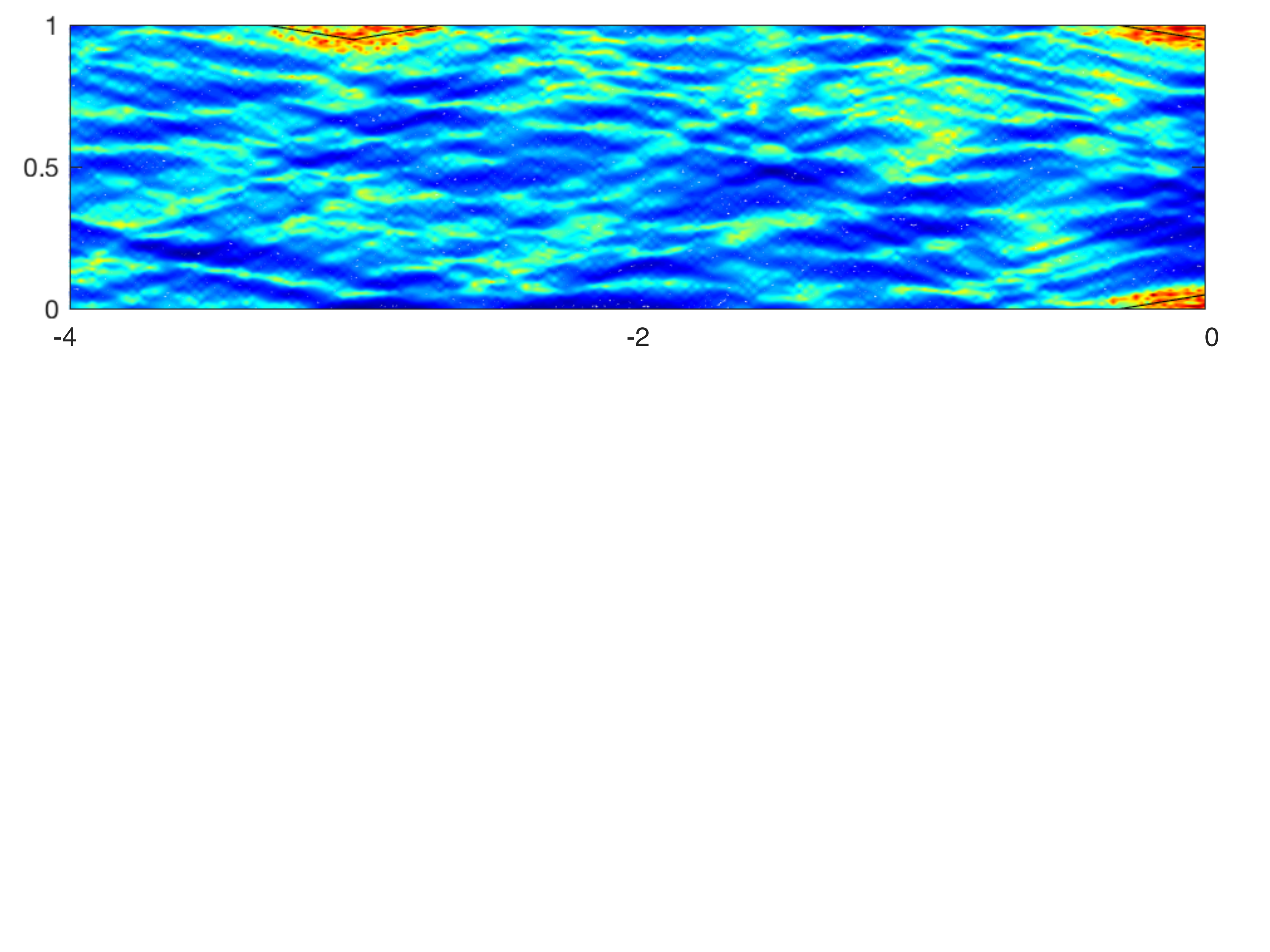}
\vspace{-4cm}
\caption{Reconstruction of wall deformations shown with a solid black
  line.  The abscissa is range and the ordinate is cross-range, scaled
  by $|\cX|$. Full aperture data 
and from top to bottom: { $J +1= 10$, $20$ and $50$}. } \label{side3_50modes}
\end{center}
\end{figure}

\begin{figure}[H]
\begin{center}
\includegraphics[scale=0.25]{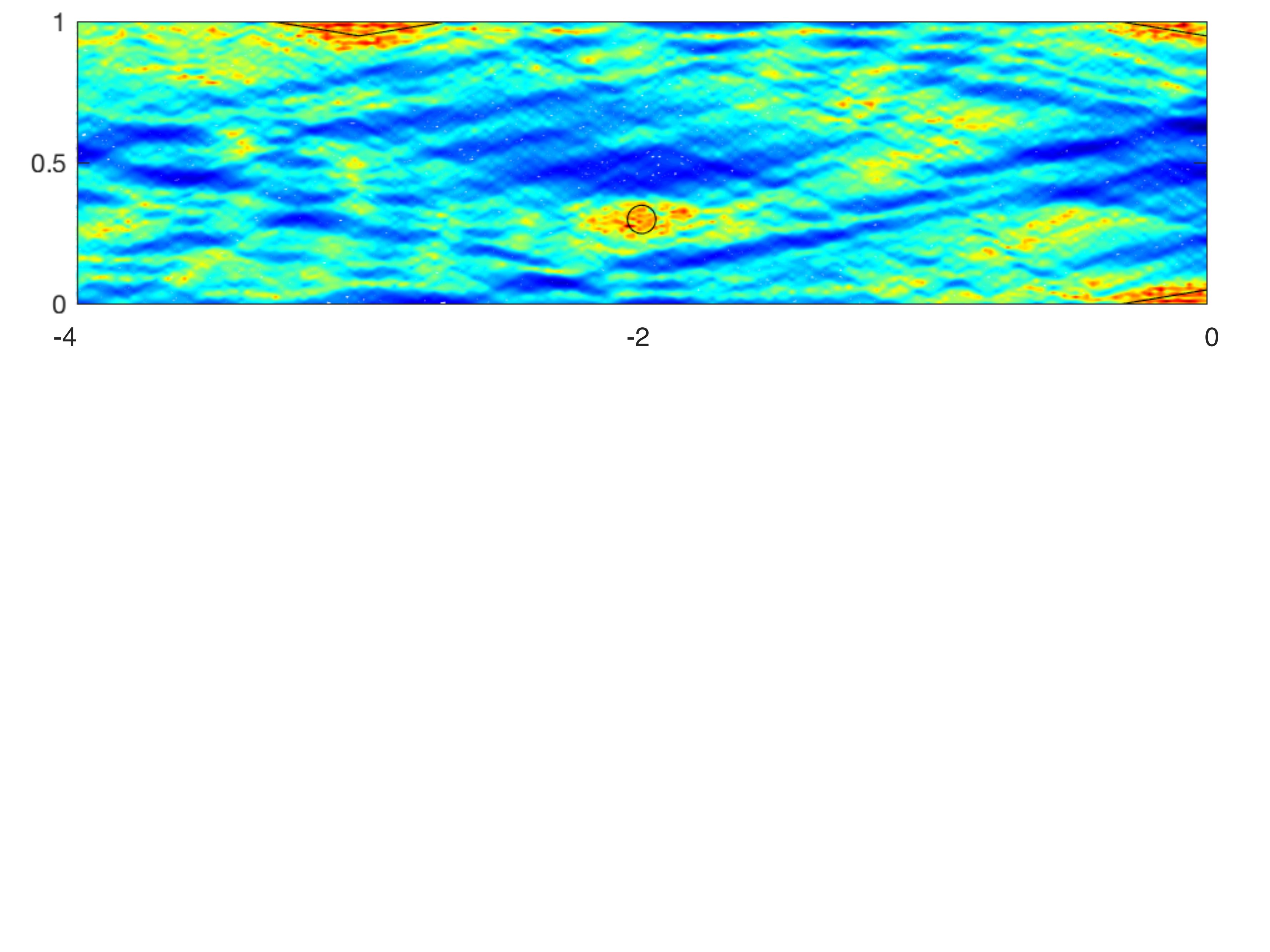}
\vspace{-4cm}
\end{center}
\begin{center}
\includegraphics[scale=0.25]{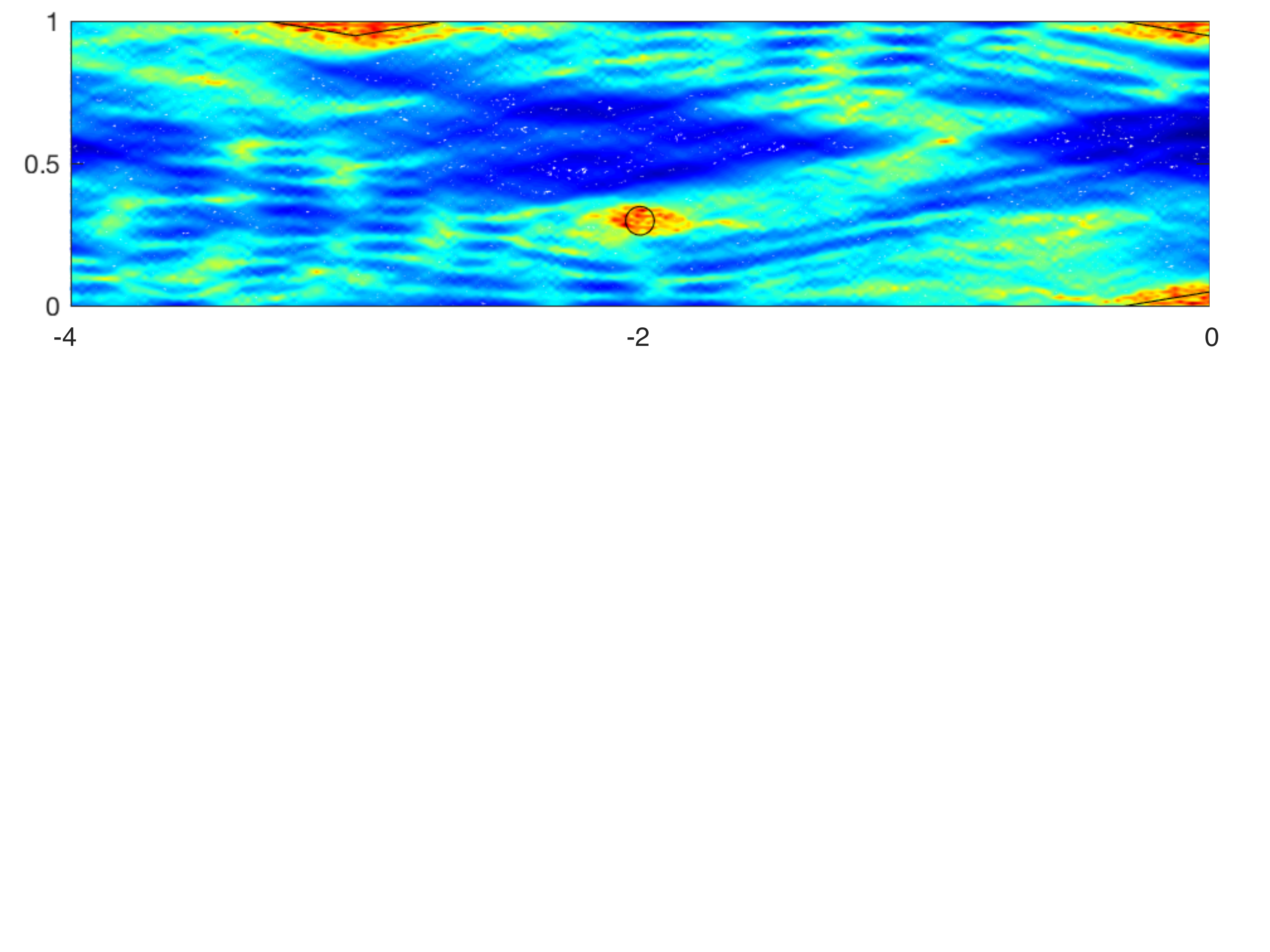}
\vspace{-4cm}
\caption{Reconstruction of wall deformations and a scatterer shown
  with a solid black line.  The abscissa is range and the ordinate is
  cross-range, scaled by $|\cX|$. Full aperture array data and {$J+1 =
  50$}. Top: sound-soft scatterer. Bottom: penetrable
  scatterer. }
\label{WallTarget}
\end{center}
\end{figure}

The first results, in Fig. \ref{side1_10modes}--\ref{WallTarget}
are obtained with a full aperture. In Fig. \ref{side1_10modes} we
show the reconstruction of wall deformations near the end of the
waveguide, for a lower frequency probing wave corresponding to { $J+1 =
10$} propagating modes. The resolution improves at higher frequencies,
as illustrated in Fig. \ref{side3_50modes}, where we show
reconstructions of wall deformations using { $J+1 = 10$, $20$ and $50$}
propagating modes.

In Fig. \ref{WallTarget} we display images in a waveguide with
wall deformations and a scatterer inside. The waveguide supports {$50$} propagating modes. 
The scatterer is impenetrable, with sound
soft boundary in the top plot, and it is penetrable in the bottom plot.

\begin{figure}[H]
\begin{center}
\includegraphics[scale=0.25]{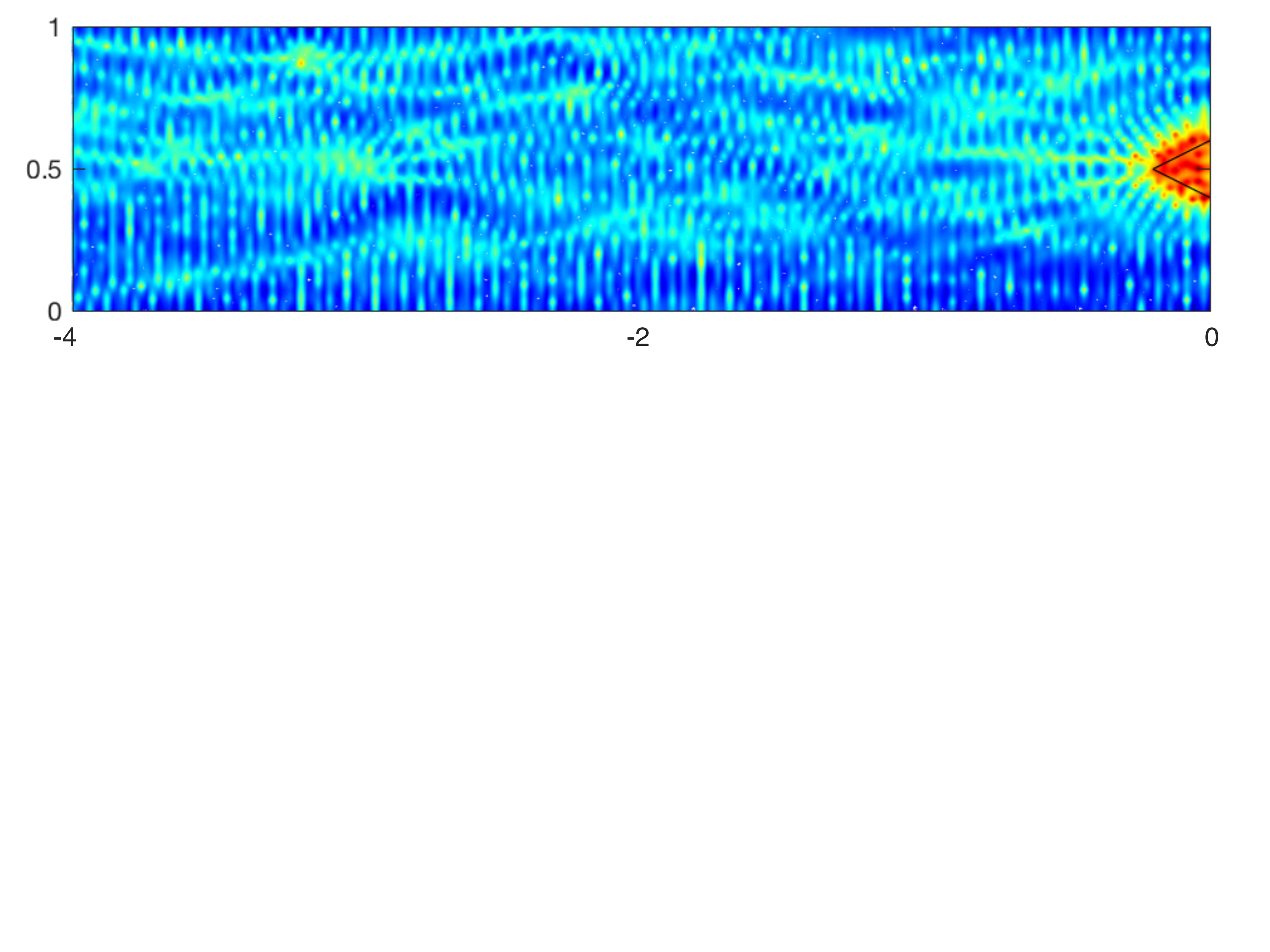}
\vspace{-4cm}
\end{center}
\begin{center}
\includegraphics[scale=0.25]{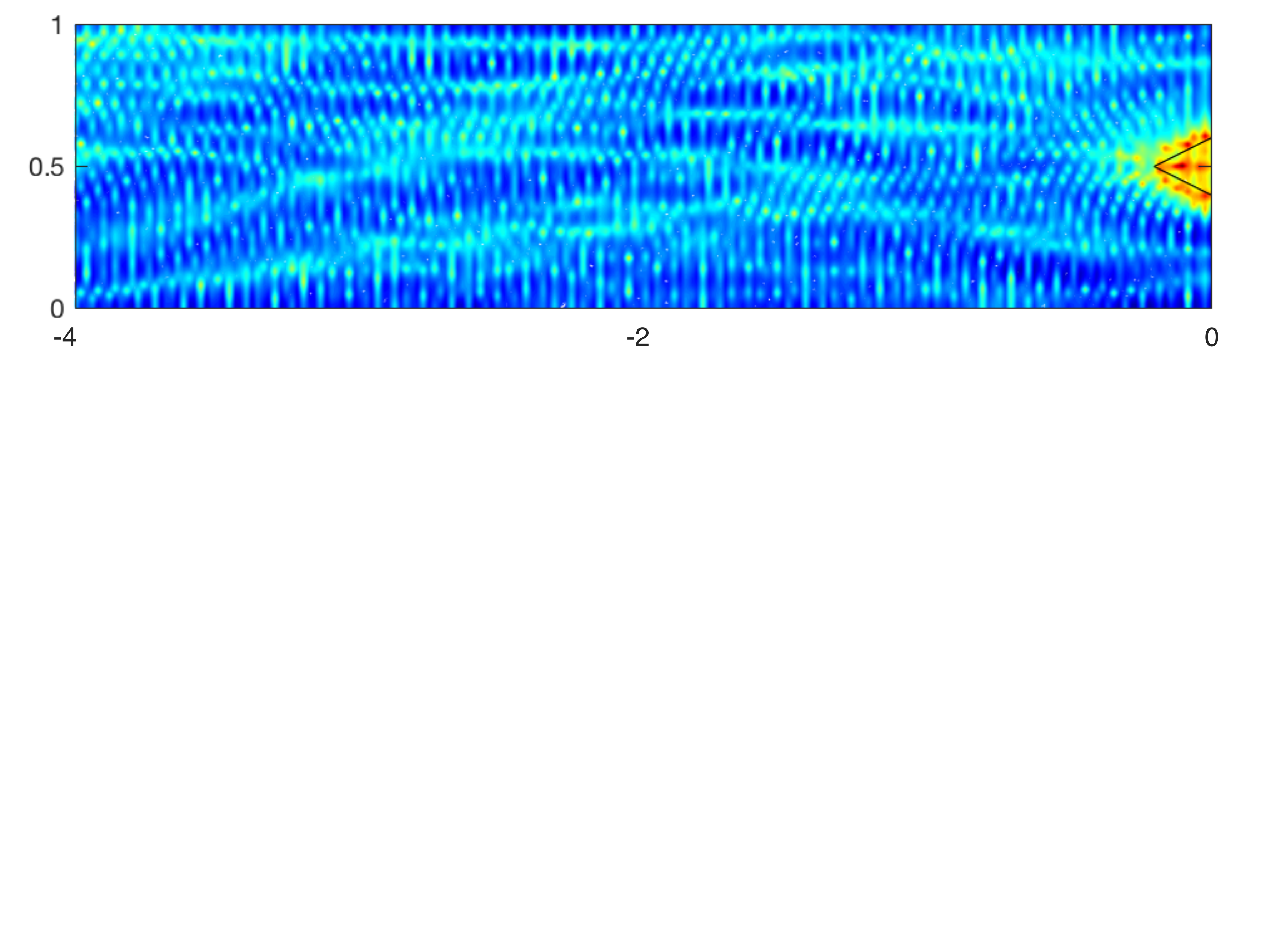}
\vspace{-4cm}
\end{center}
\begin{center}
\includegraphics[scale=0.25]{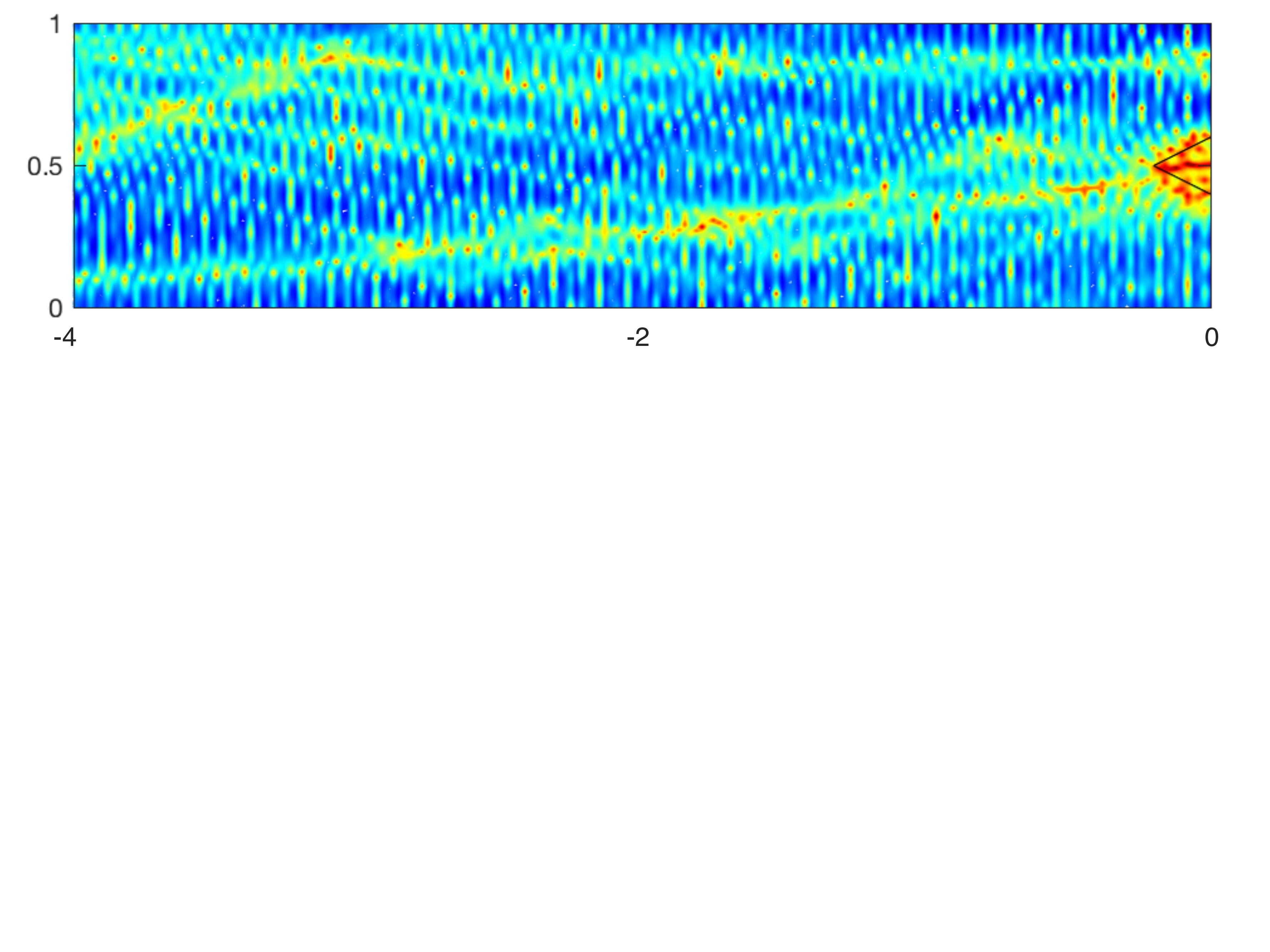}
\vspace{-4cm}
\caption{Reconstruction of wall deformations shown with a solid black
  line.  The abscissa is range and the ordinate is cross-range, scaled
  by $|\cX|$. Waveguide with $J +1= 20$. From top to bottom: $90\%$, $60\%$ and $40\%$ aperture.}
\label{PartAp}
\end{center}
\end{figure}

The effect of the aperture is illustrated in Fig. \ref{PartAp}, in
the waveguide considered in the top plot of Fig.
\ref{side1_10modes}, but this time the number of propagating modes is
increased to {$20$}. As expected, the image is better for the larger
aperture, but even when $|\cX_{\cA}|/|\cX| = 0.4$, the wall
deformation is clearly seen.

\section{Summary}
\label{sect:Summary}
We analyzed a direct approach to imaging in a waveguide with reflecting walls and  perturbed geometry. The perturbation 
consists of localized wall deformations that are unknown and are to be determined as part of the imaging. The waveguide 
may be empty or it may contain some localized, unknown  scatterers. The data  are gathered by an array of sensors 
that emits time harmonic probing waves and measures the scattered waves. Ideally, the array spans the entire cross-section of the 
waveguide, but we also consider partial aperture arrays. Starting from first principles, we established a mathematical foundation 
of the imaging algorithm. We also assessed its performance using numerical simulations in a two dimensional waveguide.

\section*{Acknowledgments}
This material is based upon research supported in part by the Air
Force Office of Scientific Research under awards FA9550-18-1-0131 and  {FA9550-17-1-0147}.

% ----------------
\appendix
\section{Proofs of Lemmas \ref{lem.1}--\ref{lem.3}}
\label{ap:A}
We analyze first in Section \ref{ap:A.1} the forward problem
\eqref{eq:1.7}--\eqref{eq:1.8} for the scattered wave field. Then we
prove the Lemmas \ref{lem.1}--\ref{lem.3} in Sections
\ref{ap:A.2}--\ref{ap:A.4}.

\subsection{Forward problem}
\label{ap:A.1}
Let us introduce the truncated waveguide 
\begin{equation}
\cW_L = \cW \cap (x_L,0) \times \cX, \qquad x_L < \xA < 0,
\label{eq:A1}
\end{equation}
between the wall at $x = 0$ and the truncation boundary
\begin{equation}
 \cL = \{x_L\} \times \cX,
\end{equation}
and show that solving the problem \eqref{eq:1.7}--\eqref{eq:1.8} in
the unbounded $\cW$ is equivalent to solving the following boundary
value problem in $\cW_L$:
\begin{align}
\big(\Delta_{\vx} + k^2 \big)u(\vx,\vx_s) &= 0, \qquad \vx \in
\cW_L, \label{eq:A2} \\ \frac{\partial u(\vx,\vx_s)}{\partial \vnu} &=
0, \qquad \vx \in { \partial \cW_L \setminus \overline{\Gamma}_o}, \label{eq:A3}
\\ \frac{\partial u(\vx,\vx_s)}{\partial \vnu} &= - \frac{\partial
  G(\vx,\vx_s)}{\partial \vnu}, \qquad \vx \in \GD, \label{eq:A4}
\\ \frac{\partial u(\vx,\vx_s)}{\partial \vnu} &= \Lambda_k
u(\vx,\vx_s), \qquad \vx \in \cL.  \label{eq:A5}
\end{align}
Here we introduced the Dirichlet to Neumann map
\begin{equation}
\Lambda_k : \widehat H^{\frac{1}{2}}(\cL)\to \widehat H^{-\frac{1}{2}}(\cL), \qquad
\Lambda_k g(\vx)\big|_{\cL} = \sum_{j=0}^\infty i \beta_j g_j
\psi_j(\bxp), \label{eq:A6}
\end{equation}
defined for all $g \in \widehat H^{\frac{1}{2}}(\cL)$, with components
\begin{equation}
g_j = \int_{\cX} d \bxp \, \psi_j(\bxp)g((x_L,\bxp)).
\label{eq:A7}
\end{equation}
The subspaces $\widehat H^{\frac{m}{2}}(\cL)$ of
$H^{\frac{m}{2}}(\cL)$ for $m = \pm 1$ correspond to functions that
satisfy Neumann boundary conditions,
\begin{equation}
\hspace{-0.1in}\widehat{H}^{\frac{m}{2}}(\cL) = \mbox{closure}\Big\{ v(\bxp) \in
\mbox{span}\{\psi_j(\bxp), ~ j \ge 0 \} ~ \mbox{s.t.}~
\sum_{j=0}^\infty (1+\lambda_j)^{\frac{m}{2}} |v_j|^2 < \infty\Big\}, 
\label{eq:A8}
\end{equation}
where 
\begin{equation}
v_j = \int_{\cX} d \bxp \, \psi_j(\bxp) v(\bxp).
\label{eq:A9}
\end{equation}
The norm in $\widehat H^{\frac{m}{2}}(\cL)$ is 
\begin{equation}
\|v\|_{\widehat H^{\frac{m}{2}}(\cL)} = \Big[ \sum_{j=0}^\infty
  (1+\lambda_j)^{\frac{m}{2}} |v_j|^2 \Big]^{\frac{1}{2}}
\label{eq:A10}
\end{equation}
and the duality pairing between $\widehat H^{-\frac{m}{2}}(\cL)$ and
$\widehat H^{\frac{m}{2}}(\cL)$ is
\begin{equation}
\left< v,w \right> = \sum_{j=0}^\infty v_j^\star w_j, \qquad \forall v
\in \widehat H^{-\frac{m}{2}}(\cL), ~ ~ \forall w
\in \widehat H^{\frac{m}{2}}(\cL),
\label{eq:A11}
\end{equation}
where the star denotes complex conjugate.

\vspace{0.05in} 
\begin{lemma}
\label{lem.A1}
The map $\Lambda_k$ is bounded for any $k$.  The map $\Lambda_i$ is
negative definite and the map $\Lambda_k - \Lambda_i$ is compact.
\end{lemma}

\vspace{0.05in} \textbf{Proof:} We have by the definition \eqref{eq:A6} that 
\begin{align*}
\|\Lambda_k g\|^2_{\widehat H^{-\frac{1}{2}}(\cL)} = \sum_{j=0}^\infty
(1+\lambda_j)^{-\frac{1}{2}} |\beta_j g_j|^2 = \sum_{j=0}^\infty
(1+\lambda_j)^{\frac{1}{2}} |g_j|^2 \frac{|\beta_j|^2}{1+ \lambda_j}
\le C \|g\|^2_{\widehat H^{\frac{1}{2}}(\cL)},
\end{align*}
 where we used definition \eqref{eq:1.6} of $\beta_j$ to obtain the bound
\[
\frac{|\beta_j|^2}{1+ \lambda_j} = \frac{|k^2-\lambda_j|}{1+\lambda_j} \le C,
\]
with constant $C > 0$ independent of $j$. This shows that $\Lambda_k$
is bounded, for any $k$.

Using the duality pairing \eqref{eq:A11}, the definition
\eqref{eq:1.6} with $k$ replaced by $i$ so that $\beta_j$ becomes $i
(1 + \lambda_j)^{\frac{1}{2}},$ and
\begin{equation}
\Lambda_i g(\vx)\big|_{\cL} = -\sum_{j=0}^\infty
(1+\lambda_j)^{\frac{1}{2}} g_j \psi_j(\bxp), \qquad \forall \, g \in
\widehat H^{\frac{1}{2}}(\cL),
\label{eq:A12}
\end{equation} 
we have for all $g \in \widehat H^{\frac{1}{2}}(\cL)$ that
\begin{align*}
\left <\Lambda_i g, g \right > = -\sum_{j=0}^\infty
(1+\lambda_j)^{\frac{1}{2}} |g_j|^2 = -\| g \|^2_{\widehat
  H^{\frac{1}{2}}(\cL)},
\end{align*}
so $\Lambda_i$ is negative definite. 

We also have from \eqref{eq:A6} and \eqref{eq:A12} that 
\begin{equation}
(\Lambda_k - \Lambda_i) g(\vx)\big|_{\cL} = \sum_{j=0}^\infty \big[i
    \beta_j + (1+\lambda_j)^{\frac{1}{2}}\big] g_j \psi_j(\bxp),
  \qquad \forall \, g \in \widehat H^{\frac{1}{2}}(\cL),
\label{eq:A13}
\end{equation} 
and we now show that in fact $(\Lambda_k - \Lambda_i) g \in \widehat
H^{\frac{1}{2}}(\cL)$. Then, the compact embedding of $\widehat
H^{\frac{1}{2}}(\cL)$ in $\widehat H^{-\frac{1}{2}}(\cL)$ gives that
$\Lambda_k - \Lambda_i$ is compact.

Indeed, we have
\begin{align}
\big\|(\Lambda_k - \Lambda_i) g\big\|^2_{\widehat
  H^{\frac{1}{2}}(\cL)} &= \sum_{j=0}^\infty (1 + \la_j)^{\frac{1}{2}}
\big|i \beta_j + (1+\la_j)^{\frac{1}{2}}\big|^2 |g_j|^2 \nonumber
\\ &= \sum_{j=0}^\infty (1 + \la_j)^{-\frac{1}{2}} |g_j|^2 \left|
\Big(i \beta_j + \sqrt{1 + \lambda_j}\Big)\sqrt{1+\lambda_j}\right|^2
\nonumber \\ &\le C \|g\|^2_{\widehat H^{-\frac{1}{2}}(\cL)},
\label{eq:A14}
\end{align}
for some positive constant $C$, because
\begin{align}
\Big|\Big( i \beta_j + \sqrt{1 + \la_j}\Big)\sqrt{1+\la_j}\Big| \le
C_1, \qquad 0 \le j \le J,
\label{eq:A16}
\end{align}
and 
\begin{align}
\Big( i \beta_j + \sqrt{1 + \la_j}\Big)\sqrt{1+\la_j} &= \Big( \sqrt{
  1+\la_j } - \sqrt{\la_j - k^2} \Big) \sqrt{1+\la_j} \nonumber \\ &=
\frac{k^2+1}{1 + \sqrt{\la_j - k^2}/\sqrt{\la_j + 1}} \nonumber \\ &\le
C_2, \qquad j > J,
\label{eq:A15}
\end{align}
where $C_1$ and $C_2$ are positive constants. Thus, 
\eqref{eq:A14} holds with $C = \max\{C_1^2,C_2^2\}$.
\endproof

\subsubsection{Connection between the scattering problems in $\cW$ and $\cW_L$}
Since problem \eqref{eq:1.7}--\eqref{eq:1.8} is stated in the infinite
domain $\cW$ and problem \eqref{eq:A2}--\eqref{eq:A5} is stated in the
truncated domain $\cW_L$, we need the following lemma to make the
connection:

\vspace{0.05in} 
\begin{lemma}
\label{lem.A2}
Consider an arbitrary $f \in \widehat H^{\frac{1}{2}}(\cL)$ with the
decomposition
\begin{equation}
f(\vx)\big|_{\cL} = \sum_{j=0}^\infty f_j \psi_j(\bxp).
\label{eq:A17}
\end{equation}  
There exists a unique solution $w \in H^1_{\rm loc}\big((-\infty, x_L)
\times \cX\big)$ of the problem
\begin{align}
\big(\Delta_{\vx} + k^2 \big)w(\vx)&= 0, \qquad \vx \in (-\infty, x_L)
\times \cX, \label{eq:A18} \\ \frac{\partial w(\vx)}{\partial \vnu} &=
0, \qquad \vx \in (-\infty, x_L) \times \partial \cX \label{eq:A19}
\\ w(\vx) &= f(\vx) , \qquad \vx \in \cL, \label{eq:A20}
\end{align}
that satisfies a radiation condition as in Definition \ref{def.1}.
\end{lemma}

\vspace{0.05in} \textbf{Proof:} From the radiation condition we know
that $w$ is an outgoing and bounded wave that has the decomposition
\begin{equation}
w(\vx) = \sum_{j=0}^\infty \gamma_j e^{-i \beta_j x} \psi_j(\bxp), \qquad 
\forall \, \vx = (x,\bxp), ~ ~ x < x_L, ~ \bxp \in \cX.
\label{eq:A21}
\end{equation}
This is a solution of \eqref{eq:A18}--\eqref{eq:A20} if 
\begin{equation}
\gamma_j = f_j e^{i \beta_j x_L}, \qquad j \ge 0,
\label{eq:A22}
\end{equation}
so the expression \eqref{eq:A21} becomes
\begin{equation}
w(\vx) = \sum_{j=0}^\infty f_j e^{-i \beta_j (x-x_L)} \psi_j(\bxp).
\label{eq:A23}
\end{equation}
Let us check that this is a function in $H^1_{\rm loc}\big((-\infty,
x_L) \times \cX\big)$. 

We have, for any $\xi < x_L$, by the orthonormality of the eigenbasis
$\{\psi_j\}_{j \ge 0}$ that
\begin{align}
\|w\|^2_{\big((-\xi, x_L) \times \cX\big)} &= \int_{\xi}^{x_L} d x
\int_{\cX} d \bxp |w(\vx)|^2 \nonumber \\ &= (x_L-\xi) \sum_{j=0}^J
|f_j|^2 + \sum_{j = J+1}^\infty |f_j|^2 \int_{\xi}^{x_L} d x \,
e^{2|\beta_j|(x-x_L)} \nonumber \\ &\le C \sum_{j=0}^\infty |f_j|^2 =
C \| f \|^2_{L^2(\cL)} \le C \|f\|_{\widehat
  H^{\frac{1}{2}}(\cL)}^2,
\label{eq:A25}
\end{align}
where $C$ is a positive constant that depends on $\xi$. Furthermore, 
using 
\begin{equation}
\nabla_{\vx} w(\vx) =  \sum_{j=0}^\infty f_j e^{-i \beta_j (x-x_L)}  
\Big(-i \beta_j\psi_j(\bxp), \nabla \psi_j(\bxp) \Big),
\label{eq:A24}
\end{equation}
the orthogonality relation 
\[
\int_{\cX} d \bxp\, \nabla \psi_j(\bxp) \cdot \nabla \psi_{j'}(\bxp) = 
\lambda_j \delta_{j,j'},
\]
and definition \eqref{eq:1.6} of the mode wavenumbers, we obtain
\begin{align}
\|\nabla_{\vx} w\|^2_{\big((-\xi, x_L) \times \cX\big)} &=
\int_{\xi}^{x_L} d x \Big[ \sum_{j=0}^J (\beta_j^2 + \lambda_j)
  |f_j|^2 + \sum_{j = J+1}^\infty (|\beta_j|^2 + \lambda_j) |f_j|^2
  e^{2|\beta_j|(x-x_L)} \Big]\nonumber \\ & \hspace{-0.5in}= (x_L-\xi)
k^2 \sum_{j=0}^J |f_j|^2 +2 \sum_{j = J+1}^\infty
(1+\lambda_j)^{\frac{1}{2}} |f_j|^2 \frac{\big[1 - e^{-2
      |\beta_j|(x_L-\xi)}\big] (\la_j -
  \frac{k^2}{2})}{\sqrt{(\la_j-k^2)(\la_j+1)}}
\nonumber\\ & \hspace{-0.5in} \le C' \|f\|_{\widehat
  H^{\frac{1}{2}}(\cL)}^2,
\label{eq:A26}
\end{align}
for another positive constant $C'$ that depends on $\xi$. The bounds
\eqref{eq:A25}--\eqref{eq:A26} and $f \in \widehat
H^{\frac{1}{2}}(\cL)$ imply that $w \in H^1_{\rm loc}\big((-\infty,
x_L) \times \cX\big)$.

It remains to prove the uniqueness of the solution. If both $w$ and
$w'$ were solutions, then $w-w'$ would also be a solution, for $f$
replaced by $0$ in \eqref{eq:A5}. Then, the estimates
\eqref{eq:A25}--\eqref{eq:A26} give that $w-w' = 0$, so the solution
is unique. \endproof

\vspace{0.05in}
\begin{theorem}
\label{thm.A1}
The scattering problem \eqref{eq:1.7}--\eqref{eq:1.8} is equivalent to 
the problem \eqref{eq:A2}--\eqref{eq:A5}.
\end{theorem}

\vspace{0.05in} \textbf{Proof:} Suppose that $u \in H^1_{\rm
  {loc}}(\cW)$ satisfies \eqref{eq:1.7}--\eqref{eq:1.8}. Then, it has
the mode expansion
\begin{equation}
u(\vx,\vx_s) = \sum_{j=0}^\infty \alpha_j \psi_j(\bxp) e^{-i \beta_j x}, \qquad 
\forall \, \vx \in (-\infty,x_L) \times \cX,
\label{eq:A27}
\end{equation}
where we suppressed the dependence of $\alpha_j$ on $\vx_s$ in the
notation. We conclude that 
\begin{equation}
u(\vx,\vx_s)\big|_{\cL} = \sum_{j=0}^\infty \alpha_j \psi_j(\bxp) e^{-i \beta_j x_L}
\label{eq:A28}
\end{equation}
is in $\widehat H^{\frac{1}{2}}(\cL)$ and using definition \eqref{eq:A6},
\begin{equation}
\Lambda_k u(\vx,\vx_s)\big|_{\cL} = -\partial_{x_L}
u(\vx,\vx_s)\big|_{\cL} = \sum_{j=0}^\infty i \beta_j \alpha_j \psi_j(\bxp)
e^{-i \beta_j x_L},
\label{eq:A29}
\end{equation}
as in \eqref{eq:A5}. Thus, $u$ satisfies \eqref{eq:A2}--\eqref{eq:A5}.

Conversely, if $u \in H^1_{\rm {loc}}(\cW_L)$ solves
\eqref{eq:A2}--\eqref{eq:A5}, we can extend it to $(-\infty,x_L)
\times \cX$ using the Dirichlet to Neumann map \eqref{eq:A6} which is
defined taking into consideration the radiation condition. \endproof

\subsubsection{Variational formulation and Fredholm alternative}
Let $v \in H^1(\cW_L)$ be arbitrary. Multiplying equation \eqref{eq:A2}
by its complex conjugate $v^\star$, integrating by parts and using the 
boundary conditions \eqref{eq:A3}--\eqref{eq:A5}, we obtain
\begin{align*}
\int_{\cW_L} d \vx \, \Big[\nabla_{\vx} u(\vx,\vx_s) \cdot
    \nabla_{\vx} v^\star(\vx) - k^2 u(\vx,\vx_s) v^\star(\vx) \Big]
  -\int_{\cL} dS_{\vx}  v^\star(\vx) \Lambda_i w(\vx) \nonumber
  \\  = - \int_{\GD} dS_{\vx}  \frac{\partial
    G(\vx,\vx_s)}{\partial \vnu} v^\star(\vx).
\end{align*}
Now let us introduce the sesquilinear forms $a(\cdot,\cdot)$ and
$h(\cdot,\cdot)$ on $H^1(\cW_L) \times H^1(\cW_L)$ and the
antilinear form $\ell(\cdot)$ on $H^1(\cW_L)$, defined by
\begin{align*}
a(w,v) &= \int_{\cW_L} d \vx \, \Big[\nabla_{\vx} w(\vx) \cdot
  \nabla_{\vx} v^\star(\vx) + w(\vx) v^\star(\vx) \Big]-\int_{\cL} dS_{\vx}
  v^\star(\vx) \Lambda_i w(\vx), \\ h(w,v) &= -(k^2+1)
\int_{\cW_L} d \vx \, w(\vx) v^\star(\vx) - \int_{\cL} dS_{\vx} 
v^\star(\vx) (\Lambda_k - \Lambda_i) w(\vx), \\ \ell(v) &= -\int_{\GD}
dS_{\vx}  \frac{\partial G(\vx,\vx_s)}{\partial \vnu} v^\star(\vx),
\qquad \forall \, w,v \in H^1(\cW_L).
\end{align*}
The variational formulation of \eqref{eq:A2}--\eqref{eq:A5} is: Find
$u(\cdot,\vx_s) \in H^1(\cW_L)$ such that
\begin{equation}
a\big(u(\cdot,\vx_s),v) + h\big(u(\cdot,\vx_s),v) = \ell(v), \qquad
\forall \, v \in H^1(\cW_L).
\label{eq:A31}
\end{equation}

From Lemma \ref{lem.A1} we know that $\Lambda_i$ is negative definite,
so it is easy to see that $a(\cdot,\cdot)$ is coercive. We also know
from Lemma \ref{lem.A1} that $\Lambda_k - \Lambda_i$ is compact, so
$h(\cdot,\cdot)$ introduces a compact perturbation of
$a(\cdot,\cdot)$.  By Fredholm's alternative, the solvability of
\eqref{eq:A31} is equivalent to the uniqueness of the
solution. Moreover, we have continuous dependence of $u$ on the
incident field at $\GD$.  

\vspace{0.05in}
\begin{theorem}
\label{thm.A2}
Let $k \in \RR$ be a positive wavenumber such that 
\begin{align}
\big(\Delta_{\vx} + k^2 \big)w(\vx) &= 0, \qquad \vx \in
\cW_L, \label{eq:A32} \\ \frac{\partial w(\vx)}{\partial \vnu}
&= 0, \qquad \vx \in { \partial \cW_L \setminus 
\overline{\cL}}, \label{eq:A33} \\ \frac{\partial
  w(\vx)}{\partial \vnu} &= \Lambda_k w(\vx), \qquad \vx
\in \cL = \{x_L\} \times \cX,  \label{eq:A34}
\end{align}
has only the trivial solution $w = 0$ in $ H^1(\cW_L)$. Then, there is a 
unique solution to \eqref{eq:A2}--\eqref{eq:A5}, and by Theorem \ref{thm.A1}
to \eqref{eq:1.7}--\eqref{eq:1.8}, and it satisfies
\begin{equation}
\| u (\cdot,\vx_s) \|_{H^1(\cW_L)} \le C_L \left\| \frac{\partial
  G(\cdot,\vx_s)}{\partial \vnu} \right\|_{H^{-\frac{1}{2}}(\GD)},
\label{eq:A35}
\end{equation}
where $C_L$ is a positive constant that depends on $x_L$.
\end{theorem}

\subsection{Proof of Lemma \ref{lem.1}}
\label{ap:A.2}
Now that we proved the solvability of the forward problem
\eqref{eq:1.7}--\eqref{eq:1.8}, we can use the definition 
of $\TDA$ in Lemma \ref{lem.1} to write 
\begin{equation}
u(\vx_r,\vx_s) = \left[\TDA \frac{\partial
  G(\cdot,\vx_s)}{\partial \vnu}\Big|_{\GD}\right](\vx_r), \qquad \vx_r \in \cA.
\label{eq:A36}
\end{equation}
Substituting in the expression \eqref{eq:1.10} of $N$ we get
\begin{equation}
N g(\vx_r) = \int_{\cA} dS_{\vx_s} \left[\TDA \frac{\partial
    G(\cdot,\vx_s)}{\partial \vnu}\Big|_{\GD}\right](\vx_r) g(\vx_s),
\qquad \forall \, g \in L^2(\cA).
\end{equation}
The integrand is smooth, so we can pull out of the integral $\TDA$ and
the normal derivative and obtain
\begin{equation}
N g(\vx_r) = \TDA \left[\partial \vnu \int_{\cA} dS_{\vx_s}
  G(\cdot,\vx_s)\Big|_{\GD}g(\vx_s)\right](\vx_r) = \TDA \TAD
g(\vx_r),
\end{equation}
where we used the definition of $\TAD$ in Lemma \ref{lem.1}. \endproof

\subsection{Proof of Lemma \ref{lem.2}}
\label{ap:A.3}
Suppose first that $\vz \in \cD$ and let $w(\vx)$ satisfy
\eqref{eq:1.13a}--\eqref{eq:1.13d} with
\[
f(\vx) = -\frac{\partial G(\vx,\vz)}{\partial \vnu}, \qquad \vx \in \GD,
\]
and the radiation condition as in Definition \ref{def.1}. Then, 
\[
v(\vx) = w(\vx) - G(\vx,\vz)
\]
solves \eqref{eq:1.13a}--\eqref{eq:1.13d} with $f = 0$. By Theorem
\ref{thm.A2}, this means that $v = 0$ so taking its trace on $\cA$ 
we get 
\[
v(\vx) = 0 = w(\vx) - G(\vx,\vz), \quad \vx \in \cA.
\]
But $w\big|_{\cA} = \TDA f$, so we have shown that
\[
\TDA f (\vx) = G(\vx,\vz), \qquad \vx \in \cA,
\]
or, equivalently, that $G(\vx,\vz)\big|_{\cA} \in \mbox{range}(\TDA)$.

Now let $\vz \notin \cD$ and suppose for a contradiction argument that
$G(\vx,\vz)\big|_{\cA}$ is in $\mbox{range}(\TDA)$. Then, there must exist 
$f \in H^{-\frac{1}{2}}(\GD)$ such that 
\[
\TDA f(\vx) = w(\vx) = G(\vx,\vz), \qquad \vx \in \cA,
\]
where $w(\vx)$ satisfies \eqref{eq:1.13a}--\eqref{eq:1.13d} and the
radiation condition. Define 
\[
v(\vx) = w(\vx) - G(\vx,\vz)
\]
and note that it satisfies 
\begin{align*}
\big(\Delta_{\vx} + k^2 \big)v(\vx) &= 0, \qquad \vx \in
(-\infty,\xA)\times \cX, \\ \frac{\partial v(\vx)}{\partial \vnu} &=
0, \qquad \vx \in (-\infty,\xA) \times \partial \cX, \\ 
v(\vx) &= 0, \quad \vx \in \cA,
\end{align*}
and the radiation condition. This problem is as in Lemma \ref{lem.A2},
with $f = 0$ and $x_L$ replaced by $\xA$. Thus, it has the unique
solution $v = 0$ in $H^1_{\rm loc}\big((-\infty,\xA) \times \cX\big)$.
By unique continuation, we can extend it to $v = 0$ in $\cW \setminus
\{\vz\}$. However, this means that $w(\vx) = G(\vx,\vz)$ which
contradicts that $w \in H^1_{\rm loc}(\cW),$ due to the singularity
of the Green's function at $\vx = \vz \in \cW$. \endproof

\subsection{Proof of Lemma \ref{lem.3}}
\label{ap:A.4}

Since $\GD$ is only part of the boundary $\partial \cW$ and $\partial
\cD$, we introduce the following Sobolev spaces on $\GD$. Suppose that
$\GD$, $ \overline{\GD} \cap (\partial \cD \backslash \GD)$ and
$\partial \cD \backslash \overline{\GD}$ are Lipschitz dissections of
the boundary $\partial \cD$. Following the notations in
\cite{mclean2000strongly}, with $\mathfrak{D}(\partial \cD) $ denoting
the space of $C^\infty(\partial \cD)$ functions with compact support,
let
\[\mathfrak{D}(\GD)=\{ \phi \in \mathfrak{D}(\partial \cD): \, 
\mbox{supp} ~ \phi \subset \Gamma_D \}.
\]
Then, we define
\begin{eqnarray*}
H^s(\GD) &=& \{ \phi|_{\GD}:\, \phi \in H^s(\partial \cD) \},
\\ \widetilde{H}^s(\GD) &=& \mbox{closure of }\,\, \mathfrak{D}(\GD)
\, \mbox{ in } \, H^s(\partial \cD),
\end{eqnarray*}
for $s = \pm\frac{1}{2}$, where the dual of {$H^s(\GD)$} is
$\widetilde{H}^{-s}(\GD)$.

Let us begin with the proof that $\TAD$ is bounded. Because
$G(\vx,\vx_s)$ is smooth for $\vx \notin \cA$, we have that
\[
v(\vx) = \int_{\cA} d S_{\vx_s}  G(\vx,\vx_s) g(\vx_s),  \qquad \forall g \in L^2(\cA),
\] 
is in $H^1\big((\xA,0) \times \cX \big)$. Moreover, 
\[
\Delta_{\vx} v(\vx) = \int_{\cA} dS_{\vx_s} \Delta_{\vx} G(\vx,\vx_s)
g(\vx_s) = - k^2 \int_{\cA} dS_{\vx_s}G(\vx,\vx_s) g(\vx_s),
\]
so we can bound
\begin{align*}
|\Delta_{\vx} v(\vx)| \le k^2 \int_{\cA} dS_{\vx_s}|G(\vx,\vx_s)
g(\vx_s)| \le C \|g\|_{L^2(\cA)},
\end{align*}
with some positive constant $C$. Here we used the Cauchy-Schwartz
inequality and that $G(\vx,\vx_s)$ is bounded for $\vx \notin
\cA$. Then, we conclude from \cite[Theorem 5.7]{mclean2000strongly} or
\cite[Lemma 4.3]{cakoni2016qualitative} that $\TAD g \in
H^{-\frac{1}{2}}(\GD)$ and its norm is bounded by the
$\|g\|_{L^2(\cA)}$.  This shows that the linear operator $\TAD$ is
bounded.

To prove that $\TAD$ has dense range in $H^{-\frac{1}{2}}(\GD)$, we
show that $h \in \widetilde H^{\frac{1}{2}}(\GD)$ must be zero if
\[
(\TAD g, h) = 0, \qquad \forall g \in L^2(\cA),
\]
where $(\cdot,\cdot)$ denotes the duality pairing. Indeed if
\begin{eqnarray*}
(\TAD g, h) = \int_{\GD} dS_{\vz} \, h^\star(\vz)
  \frac{\partial}{\partial \vnuz}\int_{\cA} dS_{\vx_s} G(\vz,\vx_s)
  g(\vx_s) =0,
\end{eqnarray*}
for all $g \in L^2(\cA)$ and $h^\star$ is the conjugate of $h$, then by the reciprocity relation
$G(\vz,\vx_s) = G(\vx_s,\vz)$ and Fubini's theorem we conclude
\begin{eqnarray*}
\int_{\GD} dS_{\vz} \, h^\star(\vz) \frac{\partial
  G(\vx_s,\vz)}{\partial \vnuz} =0, \qquad \forall \vx_s \in \cA.
\end{eqnarray*}
Let us define
\[
w(\vx) = \int_{\GD} dS_{\vz} \, h^\star(\vz) \frac{\partial
  G(\vx,\vz)}{\partial \vnuz}, \] and consider first $\vx \in
(-\infty, \xA) \times \cX$.  By Lemma \ref{lem.A2}, with $\cL$
replaced by $\cA$ and right hand side in \eqref{eq:A20} replaced by
$0$, we conclude that $w=0$ in $(-\infty, \xA) \times \cX$. Then,
unique continuation yields that
\begin{eqnarray*}
w(\vx) = \int_{\GD} dS_{\vz} \, h^\star(\vz) \frac{\partial
  G(\vx,\vz)}{\partial \vnuz} = 0, \qquad \forall \vx \in \cW.
\end{eqnarray*}
Since the Green function $G(\vx,\vz)$ has the same regularity
properties as the Green function for free space
\cite{bourgeois2008linear}, by the continuity of the double-layer
potential \cite{mclean2000strongly} we conclude that $w$ satisfies
\begin{align*}
\Delta_{\vx} w(\vx) + k^2 w(\vx) &= 0, \qquad \vx \in \cD, \\
\frac{\partial w(\vx)}{\partial \vnu} &=0, \qquad \vx \in \partial \cD,
\end{align*}
Assuming that $-k^2$ is not an eigenvalue of the Laplacian in $\cD$,
we conclude that $w=0$ in $D$.  Then, from the jump relations for
double-layer potentials (see for instance \cite{mclean2000strongly})
\begin{eqnarray*}
h^\star(\vz) =w^+(\vz) - w^-(\vz)=0, \qquad \forall \, \vz \in \GD.
\end{eqnarray*}
This concludes the proof that $\TAD$ has dense range in  $H^{-\frac{1}{2}}(\GD)$.

Now let us study the operator $\TDA$ defined in Lemma \ref{lem.1}. For
all $g \in \widetilde{H}^{\frac{1}{2}}(\GD)$ let
\[
w_g(\vx) = \int_{\GD} dS_{\vz} \, \frac{\partial G(\vx,\vz)}{\partial
  \vnuz} g(\vz), \qquad \vx \in \cW,
\]
and use the jump relations of double layer potentials to define
\[
f_g(\vx) = -\frac{\partial}{\partial \vnu} \int_{\GD} dS_{\vz} \,
\frac{\partial G(\vx,\vz)}{\partial \vnuz} g(\vz), \qquad \vx \in \GD.
\]
Since $w_g$ satisfies \eqref{eq:1.13a}--\eqref{eq:1.13d} and the
radiation condition, we can write
\[
\TDA f_g(\vx) = w_g(\vx), \qquad \vx \in \cA.
\]
To prove that $\mbox{range}(\TDA)$ is dense in $L^2(\cA)$, we show
that $h \in L^2(\cA)$, satisfying
\[
(\TDA f_g,h) = 0, \qquad \forall g\in \widetilde{H}^{\frac{1}{2}}(\GD),
\]
must be zero. Here $(\cdot,\cdot)$ denotes the $L^2(\cA)$ inner
product. Indeed, if
\begin{align*}
(\TDA f_g,h) = \int_{\cA} d S_{\vx}\, h^\star(\vx) \int_{\GD} dS_{\vz}
  \frac{\partial G(\vx,\vz)}{\partial \vnuz} g(\vz) =0, \quad \forall
  \, g \in \widetilde{H}^{\frac{1}{2}}(\GD),
\end{align*}
then, by the reciprocity relation $G(\vx,\vz) = G(\vz,\vx)$ and by
Fubini's theorem we have that
\begin{eqnarray} \label{B dense dirichlet bc}
\int_{\cA} dS_{\vx} \, h^\star(\vx) \frac{\partial
  G(\vz,\vx)}{\partial \vnuz} = 0, \qquad \forall \vz \in \GD.
\end{eqnarray}
Let 
\[
v(\vz) = \int_{\cA} dS_{\vx} \, h^\star(\vx) G(\vz,\vx), \qquad \vz
\in \cW \setminus \overline{\cA}.
\]
Since $\cA$ and $\GD$ do not intersect, we have from \eqref{B dense
  dirichlet bc} that
\[
\frac{\partial v(\vz)}{\partial \vnuz} = 0, \qquad \vz \in \GD.
\]
Furthermore, from the definition of the Green's function, 
\begin{align*}
\Delta_{\vz} v(\vz) + k^2 v(\vz) &= 0, \qquad \vz \in \cD, \\
\frac{\partial v(\vz)}{\partial \vnuz} &= 0, \qquad  \vz \in \partial \cD.
\end{align*}
Assuming that $-k^2$ is not an eigenvalue of the Lapacian in $\cD$, we
conclude that $v=0$ in $\cD$. Unique continuation yields further that
$v = 0 $ in $\cW \cap (\xA,0) \times \cX$ and from the jump relations
of the single-layer potential we get that $v = 0$ in $
H^{\frac{1}{2}}(\cA).$ Then, it follows from Lemma \ref{lem.A2} that
$v = 0 \in (-\infty,\xA) \times \cX$. The function $h$ is obtained
from the jump relations for the single layer potentials
\begin{eqnarray*}
{h}^\star(\vx) = \frac{\partial v^+(\vx)}{\partial \vnu} -
\frac{\partial v^-(\vx)}{\partial \vnu} =0, \qquad \forall \, \vx \in
\cA.
\end{eqnarray*}
This proves that $\TDA$ has dense range in $L^2(\cA)$.

Finally, from the properties of the solution of
\eqref{eq:1.13a}--\eqref{eq:1.13d} and the radiation condition we have
that
\[
\| w|_{\cA} \|_{H^{\frac{1}{2}}(\cA)} = \| \TDA f_g
\|_{{{H^{\frac{1}{2}}(\cA)}}} \le C \|f_g\|_{H^{-\frac{1}{2}}(\GD)}.
\]
The compact embedding of $H^{\frac{1}{2}}(\cA)$ in $L^2(\cA)$
gives that $\TDA$ is compact. \endproof\

%%%%%%%%%% Insert bibliography here %%%%%%%%%%%%%%

\bibliographystyle{plain} \bibliography{BCM18}

\end{document}

%% file: mathmacros.tex
\DeclareMathAlphabet{\itbf}{OML}{cmm}{b}{it}
 \DeclareMathAlphabet\mathbfcal{OMS}{cmsy}{b}{n}

\def\RR{\mathbb{R}}

\def\bv{{{\itbf v}}}
\def\bx{{{\itbf x}}}
\def\vx{\vec{\bx}}
\def\bxp{\bx^\perp}
\def\by{{{\itbf y}}}
\def\vy{\vec{\by}}
\def\byp{\by^\perp}
\def\bz{{{\itbf z}}}
\def\vz{\vec{\bz}}

\newcommand{\la}{\lambda}
\newcommand{\ep}{\varepsilon}
\newcommand{\vnu}{\vec{\boldsymbol{\nu}}_{_{\vx}}}
\newcommand{\vnuz}{\vec{\boldsymbol{\nu}}_{_{\vz}}}

\newcommand{\xA}{x_{_{\hspace{-0.03in}\cA}}}
\newcommand{\JA}{J_{_{\hspace{-0.03in}\cA}}}
\newcommand{\cA}{\mathcal{A}}
\newcommand{\cL}{\mathcal{L}}
\newcommand{\cX}{\mathfrak{X}}
\newcommand{\cW}{\mathcal{W}}
\newcommand{\cD}{\mathcal{D}}
\newcommand{\GD}{\Gamma}
\newcommand{\TDA}{T^{\GD \to \cA}}
\newcommand{\TAD}{T^{\cA \to \GD}}
\newcommand{\TDOA}{T^{\GD,\Omega \to \cA}}
\newcommand{\TADO}{T^{\cA \to \GD,\Omega}}